\documentclass{article}

\usepackage[utf8]{inputenc}
\usepackage[a4paper,top=3cm,bottom=3cm,left=3cm,right=3cm,marginparwidth=1.75cm]{geometry}
\usepackage{tikz-cd, amsmath, amssymb, amsthm, mathtools}
\usepackage{adjustbox, algorithm, algorithmic, enumitem, subcaption}
\usepackage{hyperref, natbib}
\usepackage[affil-it]{authblk}

\newcommand\phiu[2]{\varphi_{{#1}{#2}}}
\newcommand\rhov[2]{\rho_{#1 #2}}

\newcommand\inc[1]{i^{#1}}
\newcommand\pro[1]{\pi^{#1}}
\newcommand\har[1]{\phi^{#1}}

\newcommand{\sbm}[1]{{\let\amp=&\left[\begin{smallmatrix}#1\end{smallmatrix}\right]}}

\DeclareMathOperator{\Ima}{Im}

\DeclareMathOperator{\im}{im\, }
\DeclareMathOperator{\VR}{VR}
\DeclareMathOperator{\PH}{PH}
\DeclareMathOperator{\Ho}{H}
\DeclareMathOperator{\Id}{Id}
\DeclareMathOperator{\rank}{rank}
\newcommand{\subU}{Z}

\DeclareMathOperator{\PD}{{ PD}}
\DeclareMathOperator{\mtb}{{ B}}
\newcommand{\barsize}[1]{\# \PD ({#1})}
\newcommand{\matching}[1]{\mathcal{P}_{#1}}

\newcommand{\tV}{\widetilde{V}}

\newcommand{\hV}{\widehat{V}}
\newcommand{\hU}{\widehat{U}}

\newcommand{\II}{\mathcal{I}}

\newcommand{\cM}{\mathcal{M}}

\newcommand{\cA}{\mathcal{A}}
\newcommand{\cB}{\mathcal{B}}

\newcommand{\cI}{\mathcal{I}}
\newcommand{\cK}{\mathcal{K}}

\newcommand{\cS}{\mathcal{S}}

\newcommand{\tcA}{\widetilde{\cA}}

\newcommand{\hcA}{\widehat{\cA}}
\newcommand{\hcB}{\widehat{\cB}}
\newcommand{\gen}[1]{{\langle #1\rangle}}

\newcommand{\bR}{\mathbb{R}}

\newcommand{\oviot}{\overline{\iota}}

\DeclareMathOperator {\pivot} {pivot}

\newcommand{\collapseEdges}{\textbf{collapse\_edges}}
\newcommand{\flagComplex}{\textbf{flag\_complex}}
\newcommand{\computePersistence}{\textbf{compute\_persistence}}
\newcommand{\getPmMatrix}{\textbf{get\_pm\_matrix}}
\newcommand{\columnReduction}{\textbf{column\_reduction}}
\newcommand{\getMatching}{\textbf{get\_matching}}

\newcommand{\repsCollapsed}{\texttt{reps\_collapsed}}
\newcommand{\prevtime}{\texttt{prev\_time}}
\newcommand{\newtime}{\texttt{new\_time}}
\newcommand{\gammanew}{\gamma_{\text{new}}}
\newcommand{\repsL}{\texttt{reps\_L}}
\newcommand{\repsK}{\texttt{reps\_K}}
\newcommand{\imReps}{\texttt{image\_reps}}
\newcommand{\barcodeL}{\texttt{barcode\_L}}
\newcommand{\barcodeK}{\texttt{barcode\_K}}

\newcommand{\leqU}{\leq^\bullet}
\newcommand{\leqV}{\leq_\bullet}

\newcommand{\vertices}{\mathcal{V}}

\newcommand{\exend}{\unskip\null\hfill\ensuremath{\diamond}}
\graphicspath{{figures/}}
\usepackage{subcaption}

\newtheorem{theorem}{Theorem}[section]
\newtheorem{proposition}[theorem]{Proposition}
\newtheorem{lemma}[theorem]{Lemma}
\newtheorem{corollary}[theorem]{Corollary}
\newtheorem{definition}[theorem]{Definition}
\newtheorem{example}[theorem]{Example}
\newtheorem{remark}[theorem]{Remark}

\graphicspath{{figures/}}

\title{Additive Partial Matchings Induced by Persistence Maps}
\date{\vspace{-5ex}}

\author[1]{Rocio Gonzalez-Diaz\thanks{\texttt{rogodi@us.es}}}
\author[2]{Manuel Soriano-Trigueros\thanks{\texttt{msoriano4@us.es}}}
\author[1]{Álvaro Torras-Casas\thanks{\texttt{atorras@us.es}}}

\affil[1]{Dpto. Matematica Aplicada I, Universidad de Sevilla, Spain}
\affil[2]{Dpto. Ciencias de la Computación e Inteligencia Artificial, Universidad de Sevilla, Spain}

\begin{document}
	
\maketitle

\begin{abstract}
Persistent homology is a fundamental tool in Topological Data Analysis.
The associated algebraic structure is the persistence module, a sequence of vector spaces connected by linear maps.
Persistence modules admit a complete and fast-to-compute invariant known as the persistence diagram.
However, this is no longer the case for maps between persistence modules (i.e. persistence maps).
We propose a new invariant for persistence maps, consisting of a partial matching between the persistence diagrams of the domain and codomain modules.
We show that this invariant is additive with respect to the direct sum decomposition of persistence maps, is more discriminative than the image module, and is computable in cubic time.
Furthermore, we provide an implementation and demonstrate its efficiency by integrating it with edge collapse techniques for flag complexes (e.g., Vietoris–Rips complexes). As a key technical contribution, we describe how to induce a persistence map between two flag complexes that have been independently simplified via edge collapses, even when a direct simplicial map between them is no longer available.
\end{abstract}

\section{Introduction}
\label{sec:introduction}

Persistence modules have been extensively used to study persistent homology, one of the main tools of Topological Data Analysis (TDA). Informally, a (finite) persistence module is defined as a finite sequence of finite-dimensional vector spaces over a fixed field $k$, and linear maps between them, 
\[ 
\begin{tikzcd}
V_1 \arrow[r] &  V_2 \arrow[r] & \ldots \arrow[r] & V_n\,.
\end{tikzcd}
\]
This sequence arises when a dataset is studied using a filtration of complexes that evolves across a range of scale parameters.
Then, homology (with coefficients over $k$) is used to obtain the persistence module and study how the topology of these complexes evolves \citep{Computational, structure}.
Persistence modules can be decomposed up to isomorphisms as a direct sum of simpler modules known as interval modules \citep{ZC05, CrawleyBoevey2015}.
In fact, any persistence module can be characterized as a multiset of the intervals appearing in such a decomposition and their multiplicity.
This multiset is called the persistence diagram, and can be computed with cubic time complexity \citep{Computational}.
\smallskip

Maps of persistence modules (so-called persistence maps) 
appear,
for example,
when considering a point cloud with two labels \citep{yohai_alpha, chromatic}.
More concretely, a persistence map $f:V\to U$ between two persistence modules $V$ and $U$ is a finite set of linear maps $f_i:V_i\rightarrow U_i$ for $i=1,2,\ldots, n$, making the following diagram commutative
\begin{equation}
\begin{tikzcd}
U \\
V \arrow[u, "f"] 
\end{tikzcd}
\,
=
\,
\begin{tikzcd}[/tikz/column 4/.style={column sep=-0.5em}]
U_1 \arrow[r] &  U_2 \arrow[r] & \ldots \arrow[r] & U_n \\
V_1 \arrow[u, "f_1"] \arrow[r] &  V_2 \arrow[u, "f_2"] \arrow[r] & \ldots \arrow[r] & \arrow[u, "f_n"] V_n . 
\end{tikzcd}
\end{equation}

Persistence maps can be seen as ladder modules \citep{ladder}, which implies that they can be decomposed into simpler, indecomposable modules.
However,
the set of  
indecomposable modules
is wild when $n > 5$ \citep{ladder}
and the decomposition becomes difficult to interpret and compute \citep{matrix}. 
The most common and computationally efficient invariant is the image module \citep{cok},
\[
f V
\coloneqq
    f_1 V_1 \longrightarrow f_2 V_2 \longrightarrow \ldots \longrightarrow f_n V_n.
\]
It can be calculated with cubic complexity, but it is not a complete invariant and has a limited discriminative power (see Example~\ref{exa:same_image} for more details).
\smallskip

A partial bijection between persistence diagrams is called a partial matching.
In this article, our aim is to obtain an induced partial matching $\matching{f}$ from a persistence map $f: V \rightarrow U$ which is additive, determines the image module and is computable.
More precisely, our $\matching{f}$ will satisfy the following properties:
\begin{itemize}
    \item it is additive with respect to the direct sum of persistence maps (see Example~\ref{exa:limitation}),
    \item it contains more information than the image module $f V$ (in particular the image module can be obtained from it), and
    \item it can be computed in cubic time. 
\end{itemize}
However, a difficulty in defining induced partial matchings is that the result cannot be functorial~\citep{categorification}.
\smallskip

There exists another induced partial matching in the literature, $\chi_f$ \citep{BauLes2015, BauLes2020}.
This partial matching has been used for topological bootstrap \citep{ReBo23} and image segmentation \citep{faithful}, showing its relevance to TDA.
However, the following example shows a limitation to its applicability. 

\begin{example}\label{exa:limitation}
Consider the following persistence map,
\begin{equation*}
    \begin{tikzcd}[cramped]
	U \\
	V \arrow[u, "f"]
    \end{tikzcd}
    =
    \begin{tikzcd}
        k^2 \arrow[r, "\Id"] & k^2 \arrow[r] & 0 \\
        0 \arrow[r]\arrow[u] & k^2 \arrow[r, "\sbm{1 \, 0 }",
        pos=0.4]\arrow[u, "\sbm{0 \, 0 \\ 0 \, 1}", pos=0.3] & k .\arrow[u]
    \end{tikzcd}
\end{equation*}
The persistence diagram of $U$ is $\{([1,2],2) \}$, the one of $V$ is $\{([2,2],1), ([2,3],1) \}$, and $f$ decomposes as
\begin{equation*}
    \begin{tikzcd}[cramped]
        k \arrow[r, "\Id"] & k \arrow[r] & 0 \\
        0 \arrow[r]\arrow[u] & 0 \arrow[r]\arrow[u] & 0\arrow[u]
    \end{tikzcd}
        \,\oplus\,
    \begin{tikzcd}[cramped]
        0 \arrow[r] & 0 \arrow[r] & 0 \\
        0 \arrow[r]\arrow[u] & k \arrow[r, "\Id"]\arrow[u] & k \arrow[u]
    \end{tikzcd}
        \,\oplus\,
    \begin{tikzcd}[cramped]
        k \arrow[r, "\Id"] & k \arrow[r] & 0 & \\
        0 \arrow[r]\arrow[u] & k \arrow[r]\arrow[u, "\Id"] & 0 .\arrow[u] 
    \end{tikzcd} 
\end{equation*}
Then, the expected induced partial matching, 
i.e. the matching that is consistent with the decomposition of $f$, is
\[
    \emptyset \mapsto [1,2], \hspace{0.5cm}
    [2,3] \mapsto \emptyset \ \ \text{ and }\ \ 
    [2,2] \mapsto [1,2].
\]
which is the partial matching given by $\matching{f}$, as shown in Example~\ref{exa:calculation_limitation}.
However, the partial matching given by $\chi_f$ 
is:
\[
   [2, 3] \mapsto [1, 2], \; \emptyset \mapsto [1, 2], \; \text{and } [2, 2] \mapsto \emptyset,
\]
 and not the expected one. $\exend$
\end{example} 
In particular, $\chi_f$ ignores the decomposition of $f$, which results in counterintuitive pairings.
Actually, $\chi_f$ contains the same information as the image module $f V$, so it has a restricted discriminative power as an invariant of persistence maps.
Our aim in introducing $\matching{f}$ is to overcome these limitations.
Another related concept is that of induced block functions, introduced by the authors in  \citep{InducedMatchings2022}.
However, the induced block function from \citep{InducedMatchings2022} is not a partial matching and the image module cannot be recovered from it.
\smallskip 

A construction closely related to partial matchings is persistent extensions, introduced in \citep{extension}.
However, this procedure is based on the witness complex and the Functorial Dowker Theorem, and then has different properties to ours.
Another related term introduced specifically for defining algebraic Wasserstein distances is the bar-to-bar morphism \citep{algebraic_wasserstein}.
There are also invariants defined for multidimensional persistence that can be applied to persistence maps.
We would like to highlight the so-called interval approximations \citep{asashiba_mobius, asashiba_approximation, hiraoka_refinement} and the closely related generalized rank invariant \citep{kim_generalized, clause_generalized_2} and ladder invariant \citep{jacquard_topics_2024}.
However, there is no immediate relation between these invariants and ours, and a comprehensive analysis is beyond the scope of this paper.
\smallskip

From a computational standpoint, we not only provide an algorithm to calculate the partial matching, but also explain how to efficiently calculate persistence maps from flag complexes—a common situation in practice. 
Flag complexes are simplicial complexes whose n-simplices are determined by the n-cliques of their 1-skeleton.
Our interest in studying these complexes is that, in TDA, persistence modules often arise from computing the homology of a filtration of a flag complex. 
Furthermore, there exists a technique known as edge collapse that greatly simplifies the calculation of these persistence modules \citep{edge_collapse, swap}.
\smallskip

When studying the relation between a complex and a subcomplex, we obtain a persistence map \citep{yohai_alpha, chromatic}.  
However, a difficulty that comes up when working with flag complexes is that a complex and its subcomplex may lose the inclusion relation after being independently simplified via edge collapses, making it non-trivial to construct a persistence map afterward. 
In this paper, we describe how to induce a persistence map between two flag complexes after simplifying them using edge collapses. 
This allows us to greatly reduce the computational cost of the persistence map for these specific cases.
\smallskip

We have organized the article as follows: we provide the necessary background in Section~\ref{sec:background}, and define $\matching{f}$ and provide an algorithm to calculate it in Section~\ref{sec:induced}.
We prove that $\matching{f}$ is additive with respect to the direct sum in Section~\ref{sec:properties}, and that $f V$ can be obtained directly from it in Section~\ref{sec:image}.
In Section~\ref{sec:collapses}, we describe how to induce a persistence map between two flag complexes that have been independently simplified via edge collapses.
We have also included computational experiments and examples in Section~\ref{sec:experiments}.
Lastly, a full comparison between $\matching{f}$ and $\chi_f$ can be found in \ref{app:BL}, and between $\matching{f}$ and the induced block functions in \ref{app:BF}.
\smallskip

We would also like to point out that this article extends the proceedings publication \citep{ISSAC2025}. 
In particular, Sections~\ref{sec:collapses}, \ref{sec:experiments}, \ref{app:BL}, and \ref{app:BF} are completely new.

\section{Background}
\label{sec:background}
In this section, we schematically review  the necessary background to understand the article. For a more general introduction, we refer to \citep{structure}.

\subsection{Persistence modules}
\label{sec:persistence-modules}

In the following, we assume that all vector spaces are defined over a fixed field $k$. 
For every natural number $n$,
the expression $[n]$ 
denotes the set $\{1, 2,\ldots, n\}$.
A \emph{persistence module} (p.m.) $V$ indexed by $[n]$ consists of
a finite set of vector spaces $V_p$ for $p \in [n]$, and a set of linear maps 
$\rhov{p}{q} : V_p \stackrel{}{\rightarrow} V_q$
for $p,q\in [n]$ and $p \leq q$, such that $\rhov{q}{l}\rhov{p}{q} = \rhov{p}{l}$ if $p \leq q \leq l$ and $\rhov{p}{p}$ is the identity map.
The maps $\rhov{p}{q}$ are known as the \emph{structure maps} of $V$.
To simplify notation, we add the trivial vector spaces $V_0 = V_{n+1} = 0$ to each p.m., as well as the respective zero structure maps $\rho_{0p}$ and $\rho_{p(n+1)}$ for $0\leq p\leq n+1$.
During the rest of the article, we consider the field $k$ and the value $n$ to be fixed.
\smallskip

Given $a,b \in [n]$, we write as $I=[a,b]$ the interval set
$\{a, a+1, \dots,b\}$, 
and as $\II(n)$ the set of all the intervals that are subsets of $[n]$.
The \emph{interval module} $k_I$ is composed of the vector spaces $k_{It}$ for $t \in [n]$ where $k_{It}= k$ if $t \in I$ and $k_{It}= 0$ otherwise,
and by the identity maps as structure maps whenever possible, or the zero map otherwise.
\[
    k_{[a,b]} := 
     \ldots \overset{0}{\longrightarrow} 0 \overset{0}{\longrightarrow}
    \underbrace{ k \overset{\Id}{\longrightarrow}  \ldots \overset{\Id}{\longrightarrow} k}_{[a , b]}   \overset{0}{\longrightarrow} 0
    \overset{0}{\longrightarrow}  \ldots .
\]

In the following, $V$ and $U$ will denote two persistence modules and $\rho$ and $\varphi$ their structure maps.
The \emph{direct sum} of $V$ and $U$, $V \oplus U$, is defined using the vector spaces $V_p \oplus U_{p}$ and the structure maps $\rhov{p}{q} \oplus \phiu{p}{q}$. 
We say that $V$ is a \emph{submodule} of $U$ if $V_p \subseteq U_p$ for all $p \in [n]$ and $\rhov{p}{q}  = \phiu{p}{q}|_{V_p}$ for any pair $p \leq q$. 
If $V$ is a submodule of $U$, the \emph{quotient} $U / V$ is the p.m. whose vector spaces are $U_i / V_i$ for all $i\in [n]$ and whose structure maps are induced by the structure maps of $U$.
\smallskip

As anticipated in Section~\ref{sec:introduction}, a \emph{persistence map} $f\colon V\to U$ is a set of linear maps $\{ f_p \}_{p \in [n]}$, such that the following diagram commutes,
\begin{equation}\label{eq:persistence_map}
    \begin{tikzcd}[/tikz/column 6/.style={column sep=-0.5em}]
    0  \arrow[r] &U_1 \arrow[r, "\phiu{1}{2}"] &  U_2 \arrow[r, "\phiu{2}{3}"] & \ldots \arrow[r, "\phiu{(n-1)}{n}", pos=0.4] & U_n \arrow[r]&0 \\
    0 \arrow[r] 	&V_1 \arrow[u, "f_1"] \arrow[r, "\rhov{1}{2}"] &  V_2 \arrow[u, "f_2"] \arrow[r, "\rhov{2}{3}"] & \ldots \arrow[r, "\rhov{(n-1)}{n}", pos=0.4] & \arrow[u, "f_n"] V_n\arrow[r]& 0. 
    \end{tikzcd}
\end{equation}
\smallskip

The persistence map $f$ is an \emph{isomorphism}\-/\-\emph{surjection}\-/\-\emph{injection} when $f_p$ is an isomorphism\-/\-surjection\-/\-injection of vector spaces for all $p \in [n]$.
We define the direct sum of two persistence maps as the pointwise direct sum of their linear maps.
We say that $V$ (or a persistence map $f$) is indecomposable if $V \cong V_1 \oplus V_2$ implies $V_1=0$ or $V_2=0$ (if $f \cong f_1 \oplus f_2$ implies $f_1 = 0$ or $f_2 = 0$).
\smallskip

Note that if $g \colon k_{[a,b]} \rightarrow k_{[c,d]}$ is a persistence map, then $g$ is necessarily null unless $c \leq a \leq d \leq b$.
When these inequalities are satisfied, we say that $[c,d] \leq [a,b]$.
From now on, we omit the subindices of the linear map when they are clear from the context. 
In particular, we may write $f V_p$ instead of $f_p V_p$.
In addition, for every $f$, its image module \citep{cok} is the p.m. given by
\[
f V \coloneqq
\begin{tikzcd}[column sep=3em, /tikz/column 3/.style={column sep=4.2em}]
    f V_1\; \arrow[r, "\phiu{1}{2}\mid_{f V_1}\;\;\text{ }"] & 
    \;f V_2\; \arrow[r, "\phiu{2}{3}\mid_{f V_2}\;\;\text{ }"] & 
    \;\cdots\;  \arrow[r, "\phiu{(n-1)}{n}\mid_{f V_{n-1}}\;\;\text{ }"] &
    \; f V_n.
\end{tikzcd}
\]
\begin{example}\label{ex:complex_indecomposable}
The following persistence map is
indecomposable and surjective:
 \begin{equation*}		
    \begin{tikzcd}
	U \\
	V \arrow[u, "f"]
    \end{tikzcd}
    =
    \begin{tikzcd}[/tikz/column 4/.style={column sep=-0.5em}] 
        k \arrow[r, "\sbm{1\\0}"] &  
        k^2 \arrow[r,"\sbm{1\;0}", pos=0.4] & k \arrow[r] & 
        0
        \\
	k \arrow[r, "\sbm{1 \\ 0}"]  \arrow[u, swap, "\Id"] &
        k^2 \arrow[r,swap, "\Id"]  \arrow[u, swap,"\sbm{1\;1\\ 0\;1}", pos=0.3] &
        k^2  \arrow[u, swap,"\sbm{1\;1}", pos=0.4] \arrow[r,"\sbm{1\;0}", swap, pos=0.4] &
        k, \arrow[u]
    \end{tikzcd}
\end{equation*}
 and its image module
\[
fV = 
    \begin{tikzcd}[/tikz/column 4/.style={column sep=-0.5em}]
        k \arrow[r, "\sbm{1\\0}"] &  
        k^2 \arrow[r,"\sbm{1\;1}", pos=0.4] & k \arrow[r] & 
        0
    \end{tikzcd}
\]
is isomorphic to $U$. $\exend$
\end{example}

\subsection{Persistence diagrams and partial matchings}
In order to represent persistence modules, we need the concept of \emph{multiset}. 
A \emph{multiset} is a set where elements may appear more than once.
Formally, it is defined as $\mathbf{S} = \{ (s, m_s) \}_{ s \in S} \subset S \times \mathbb{Z}_{>0}$ where $S$ is a finite set and $m_s$ represents the number of copies of $s$ in $S$, i.e. its multiplicity.
We also define the \emph{size} of $\mathbf{S}$ as
\[
 \# \mathbf{S} = \sum_{s\in S} m_s.
\]

The following theorem highlights the central role of multisets.
It states that persistence modules can be expressed uniquely (up to isomorphism) as a direct sum of interval modules.

\begin{theorem}[\cite{ZC05, CrawleyBoevey2015}]\label{the:decomposition}
    If $V$ is a persistence module,
    there exists a set of intervals $\II(V)\subset \II(n)$ (for some $n> 0$) and a multiset $\{(I,m_I^V)\} \subset \II(V) \times \mathbb{Z}_{>0}$ such that 
    \[
        V \cong \bigoplus_{I\in \II(V)
        } \left( \bigoplus_{1}^{m_I^V} k_I \right).
    \]
\end{theorem}

The multiset given by Theorem~\ref{the:decomposition} is called the \emph{persistence diagram} (PD) of $V$, and is denoted as $\PD(V)$.

\begin{example}\label{exa:diagram}
    The persistence module $V \coloneq 0 \longrightarrow k^2 \overset{\sbm{1 \, 0 }}{\longrightarrow} k$
    is isomorphic to the following direct sum of 
    indecomposable modules:
    \begin{equation*}
    \Big(
        \begin{tikzcd}[cramped]
            0 \arrow[r] & k \arrow[r] & 0
        \end{tikzcd}
    \Big)
        \quad
            \oplus
        \quad
    \Big(
        \begin{tikzcd}[cramped]
            0 \arrow[r]& k \arrow[r, "\Id"] & k
        \end{tikzcd}
    \Big).
    \end{equation*}
    Hence, $\PD(V) = \{([2,2],1), ([2,3],1)\}$. \exend
\end{example}

Relations between persistence diagrams are given by partial bijections, known as partial matchings, see~\citep{BauLes2015}.
Here we give an alternative definition, which is more flexible and will allow us to define the partial matching independently of a choice of basis. 
\begin{definition}\label{def:partial_matching}
    A \emph{partial matching} between $V$ and $U$ is a function $\matching{} \colon\II(n) \times \II(n) \rightarrow \mathbb{Z}_{\geq 0}$ such that $\matching{}(I, J) = 0$ if $I \not\in \II(V)$ or $J \not\in \II(U)$ and
    \begin{align}
    \label{eq:ine1}
    \sum_{J \in \II(U)} \matching{}(I, J) &\leq m_I^V, \\
    \label{eq:ine2}
    \sum_{I \in \II(V)
    } \matching{}(I, J) &\leq m_J^U.
    \end{align}
\end{definition}
Note that $\matching{}(I, J)$ can be interpreted as the number of copies of intervals $I \in \II(V)$ which are matched to $J \in \II(U)$, giving the connection with the partial bijection definition.
Inequalities~(\ref{eq:ine1}) and~(\ref{eq:ine2}) guarantee that an interval is never matched exceeding its multiplicity.

\subsection{Persistence bases}\label{sec:basis}
In order to represent persistence maps in matrix form, we need bases for persistence modules.
First, we define an indexing for $\PD(V)$ as $\iota \colon [\barsize{V}] \to \II(V)$ where the number of indices $\ell$ such that $\iota(\ell) = I$ equals $m_I^V$.
If $I = [a,b]$ and $\iota(\ell) = I$ we write $[a_\ell, b_\ell]$ for convenience.
In the following, we assume that the persistence diagrams are indexed by such an indexing function $\iota$ which we omit in the notation.
A \emph{persistence basis} \citep{distributed, structure}  for $V$ is an isomorphism
\[
    \alpha:
    \bigoplus_{\ell\in [\barsize{V}]}\, k_{[a_\ell, b_\ell]} \rightarrow V .
\]
Such an isomorphism always exists by Theorem~\ref{the:decomposition}.
The \emph{persistence generator} $\alpha_\ell:k_{[a_\ell, b_\ell]} \rightarrow V$ is defined as the restriction of the persistence map $\alpha$ to $k_{[a_\ell, b_\ell]}$ for $\ell \in [\barsize{V}]$.
When we write $\alpha_\ell \sim {[a_\ell, b_\ell]}$, we mean that $ k_{{[a_\ell, b_\ell]}}$ is the domain of $\alpha_\ell$.
We also specify a persistence basis $\alpha$ by its set of persistence generators $\cA=\{\alpha_\ell
\}_{\ell \in [\barsize{V}]}$.
Note that $\# \cA = \barsize{V}$.

\begin{example}\label{ex:bases_alpha}
Consider $V$ from
Example~\ref{exa:diagram}, and recall that it is isomorphic to $k_{[2,2]}\oplus k_{[2,3]}$. 
Then $\alpha:k_{[2,2]}\oplus k_{[2,3]} \rightarrow V$
is a persistence basis for $V$.
The persistence generators 
$\alpha_1$ and $\alpha_2$ are given by the following commutative diagrams:
\\
 \noindent  
 \adjustbox{scale=.9,center}{
 $
    \begin{tikzcd}
	V \\
	k_{[2,2]} \arrow[u, "\alpha_1"]
    \end{tikzcd}  
      = 
      \hspace{0.05cm}
    \begin{tikzcd} 
        0 \arrow[r] &  
        k^2 \arrow[r,"\sbm{1\;0}", pos=0.4] & 
        k 
        \\
    	0 \arrow[u]\arrow[r] &
        k  \arrow[u, swap,"\sbm{0\\1}"] \arrow[r] & 
        0, \arrow[u] 
    \end{tikzcd}
        $
\hspace{0.3cm}
    $
    \begin{tikzcd}[/tikz/column 3/.style={column sep=-0.5em}]
	V \\
	k_{[2,3]} \arrow[u, "\alpha_2"]
    \end{tikzcd}
    = 
    \hspace{0.05cm}
    \begin{tikzcd}[/tikz/column 3/.style={column sep=-0.5em}]
        0 \arrow[r] &  
        k^2 \arrow[r,"\sbm{1\;0}", pos=0.4] & 
        k 
        \\
    	0 \arrow[u]\arrow[r] &
        k  \arrow[u, swap,"\sbm{1\\0}"] \arrow[r,swap,"\Id"] & 
        k \arrow[u, "\Id"] & \exend
    \end{tikzcd}
    $
    }
\end{example}

\label{def:persistence-basis}
Given a subset $\cS = \{ \alpha_i\}_{i \in \Lambda}$ of a persistence basis
$\cA$ with indices $\Lambda \subseteq [\# \cA]$,
we define the \emph{span}  of $\cS$, denoted by $\gen{\cS}$, as the image of the sum of persistence generators in
$\cS$, that is  
\[
\gen{\cS}=\Ima\left(\bigoplus_{i \in \Lambda} \alpha_i\colon\bigoplus_{i \in \Lambda}\, k_{[a_i, b_i]} \rightarrow V \right).
\]

For $t \in [n]$, we define $\cS_t \coloneqq \{\alpha_{i t}^1 \mid i \in \Lambda \text{ and }  t\in [a_i, b_i]\}$ where $\alpha_{i t}^1$ 
stands for $\alpha_{i t}(1_k)$, denoting $1_k$ as the unit of $k$. 
In particular, $ V_t \cong \gen{\cA_t} \cong \gen{\cA}_t.$

\subsection{Persistent homology and flag complexes}\label{sec:flag-complexes}
In computational topology, persistence modules usually arise when calculating the homology of a sequence of simplicial complexes \citep{Computational}.
One of the most common types of simplicial complexes is the \emph{flag} (or \emph{clique}) complex, which can be completely described by its 1-skeleton.
Since their use in applications is widely spread, we study them for our computations in Sections~\ref{sec:collapses} and \ref{sec:experiments}. However, we note that the induced matching introduced in this work is defined algebraically and does not depend on a particular choice of filtered complexes.

\paragraph{Flag complexes} 
Given an ordered set of vertices, $\vertices = \{v_i\}_{i = 0, \ldots, m}$, a \emph{simplicial complex}  $K$ is a family of ordered subsets of $\vertices$, such that if $\tau \subset \sigma$ and $\sigma \in K$, then $\tau \in K$.
Every set in $K$ is denoted a \emph{simplex}, and if $\tau \subset \sigma$ we say that $\tau$ is a \emph{face} of $\sigma$ and $\sigma$ a \emph{coface} of $\tau$.
The \emph{dimension} of $\sigma$ is given by $\# \sigma - 1$.
When $\sigma$ has dimension $d$, we say it is a $d$-simplex.
We denote the \emph{$1$-skeleton} of $K$ as the graph formed by the $0$-simplices (vertices) and $1$-simplices (edges) of $K$.
A simplicial complex $K$ is a \emph{flag complex} when $\sigma \in K$ if and only if the vertices of $\sigma$ form a clique in the $1$-skeleton of $K$.
Hence, flag complexes are completely determined by their $1$-skeleton.

\paragraph{Homology} 
We define an \emph{oriented simplex} as a simplex with an orientation, i.e. an order defined on its vertices.
Note that since the vertices of the simplicial complex are ordered, they induce a canonical order on each simplex.
We use brackets to denote an oriented simplex,
for example, $[v_0, \ldots, v_d]$ represents the simplex $\{v_0, \ldots, v_d\}$ with the canonical orientation.
Given a simplicial complex $K$, we consider the $k$-vector space generated by all its oriented simplices of dimension $d$ (i.e. all its simplices, with all its possible orientations) and denote it as ${OC}_d(K)$.
We also consider the following equivalence relation: given two orientations of the same simplex, $\sigma, \tau$, we say that $\sigma \sim \tau$ if there exists an even permutation between their orders, and $\sigma \sim -\tau$ otherwise.
We denote by $C_d(K)$ 
the $k$-vector space obtained from the quotient ${OC}_d(K) / \sim$. 
Every vector in $C_d(K)$ is called a \emph{$d$-chain}.
For each dimension $d > 0$, we define the boundary operator $\partial_d \colon C_d(K) \to C_{d-1}(K)$ as follows:
given a simplex $[v_0,\cdots, v_d]\in C_d(K)$, we set
\[
	\partial_d([v_0, v_1, \ldots, v_d]) = \sum_{i=0,\ldots,d} (-1)^{i}[v_0, \ldots, \widehat{v_i}, \ldots, v_d]\ ,
\] 
where $\widehat{v_i}$ means we have removed $v_i$ from the simplex. Then, this definition of $\partial_d$ is extended linearly to any $d$-chain.
It can be checked that $\partial_d \circ \partial_{d+1} = 0$.
Such a pair $(C_d(K), \partial_d)_{d\geq 0}$ is known as a \emph{chain complex} (where we set $\partial_0$ to be the trivial map $0\colon C_0(K)\rightarrow 0$).
For ease, we often omit writing subindices of a chain complex $(C(K), \partial)$.
The $d$-homology of $K$, $\Ho_d(K)$, is the vector space $\ker \partial_d \, / \im \partial_{d+1}$.
In addition, recall that homology is a functor, so it induces a linear map  $\Ho_d(K) \to \Ho_d(L)$ whenever we have an inclusion $K \subset L$.

\paragraph{Persistent homology}
When we have a sequence of simplicial complexes, we can apply the homology functor to obtain a persistence module
\[
	K_1 \subset K_2 \subset \ldots \subset K_n
	\quad \Longrightarrow \quad
	\Ho_d(K_1) \to \Ho_d(K_2) \to \ldots \to \Ho_d(K_n).
\]
This process is called \emph{persistent homology} and was first introduced in \citep{diagrams}.
In Topological Data Analysis, it is common to obtain such a sequence of complexes using the Vietoris-Rips filtration \citep{Computational}.
It is calculated from a finite metric space $(X, \textbf{d})$ and a sequence of values $r_i \in \bR$.
For each $x \in X$, we construct a graph with the vertices of $X$ and edges ${x,y}$ if $\textbf{d}(x,y) \leq r_i$.
The Vietoris-Rips complex $\VR(r_i)$ is the flag complex of such graph, and the Vietoris-Rips filtration is obtained from the fact that $r_i < r_{i+1}$ implies $\VR(r_i) \subset \VR(r_{i+1})$.

\paragraph{Edge collapses}
Let $K$ be a flag complex and $e = \{u,v\}$ a 1-simplex in $K$.
Then, we say that $v' \neq u,v$ dominates $e$ if $\{v', u\}$ and $\{v',v\}$ are in $K$ and for every $w$ such that $\{u,w\}$ and $\{v,w\}$ are in $K$ then $\{ v', w\}$ is also in $K$.
In such case, removing $e$ and all its cofaces from $K$ generates a new flag complex $K'$ such that $\Ho(K) \cong \Ho(K')$ \citep{edge_collapse}.
This process is known as \emph{edge collapse} and we denote it by $K \setminus e$.
In \citep{swap}, some operations using edge collapses were introduced to simplify flag complexes and allow fast calculations of persistent homology.
The most important for this article are:
\begin{itemize}
	\item Shifting. Let $e$ be a dominated edge in $K$ inserted at $K_{i}$.
	Then, the insertion of $e$ can be shifted by one grading to $K_{i+1}$ without changing the persistence module induced by homology.
	In other words, the filtrations $K_1 \subset \ldots \subset K_{i} \subset K_{i+1} \subset \ldots$ and $K_1 \subset \ldots \subset K_{i} \setminus e \subset K_{i+1} \subset \ldots$ have the same persistent homology.
	\item Trimming. If a dominated edge is inserted  at the last step of the filtration, then it can be removed without altering the persistence module induced by homology.
\end{itemize}  

\subsection{Chain contractions}\label{sec:chain_contraction}

A chain contraction is an operation that allows us to substitute a chain complex with a smaller one with the same homology \citep{reduction}.
As we prove in Section~\ref{sec:collapses}, every edge collapse in a flag complex induces a chain contraction in the associated chain complexes.
\smallskip

Next, we review the concept of chain contraction. 
For this, recall that a \emph{chain map} $f\colon C\rightarrow C'$
(of degree $m \geq 0$)
between two chain complexes $(C,\partial)$ and $(C',\partial')$ is a sequence of linear maps $f_d\colon C_d\rightarrow C'_d$ (for all $d\geq 0$) such that $f_{d-1}\partial_d = \partial_{m+d}f_d$ for all $d>0$.
If we do not specify the degree of a chain map it is understood that its degree is $0$.

\begin{definition}\label{def:chain_contraction}
	Given the simplicial complexes $K$ and $L$, a \emph{chain contraction} is a triplet $(\pro{}, \inc{}, \har{})$  consisting of a pair of chain maps $\inc{} \colon C(L) \to C(K)$ and $\pro{} \colon C(K) \to C(L)$ together with a degree-$1$ chain map $\har{} \colon C(K) \to C(K)$ such that
	\[
		\pro{}\,\inc{} = \Id, \qquad \inc{} \, \pro{} = \Id + \partial \har{} + \har{} \partial.
	\]
	and
	\[
		\pro{} \, \har{} = 0, \quad \har{} \, \inc{} = 0.
	\]
\end{definition}

It is a direct consequence of the definition that, given a chain contraction, the maps $\inc{}_* \colon \Ho(L) \to \Ho(K)$ and $\pro{}_* \colon  \Ho(K) \to \Ho(L)$ which are respectively induced by $\inc{} \colon L\hookrightarrow K$ and $\pro{} \colon K\hookrightarrow L$, are isomorphisms, see \citep{reduction}. 
Chain contractions can also be composed as follows \citep{reduction}[Section 12].
If $(\pro{1}, \inc{1}, \har{1})$ is a chain contraction between $K$ and $L$, and  $(\pro{2}, \inc{2}, \har{2})$ a chain contraction between $L$ and $P$, then $(\pro{2} \, \pro{1}, \inc{1} \, \inc{2}, \har{1} + \inc{1} \, \har{2} \, \pro{1})$ is a chain contraction between $K$ and $P$.

\subsection{Lemmas about vector spaces}
We recall some minor lemmas about vector spaces that will become necessary throughout the text.

\begin{lemma}\label{lem:commute_intersection}
    Let $C,D,E$ be subspaces of a vector space $A$. If $D\subseteq E$ then 
    \[
    (C+D)\cap E = C\cap E + D.
    \]
\end{lemma}

\begin{proof}
    First, notice we have the inclusion $C \cap E + D = C\cap E + D\cap E \subseteq (C+D)\cap E$. On the other hand, let $w \in (C+D)\cap E$. Then, there exist $c \in C$ and $d \in D$ such that $w=c+d$ and so $c = w-d \in E$, since $D\subseteq E$. Thus, $w = c + d \in C\cap E + D$ as claimed.
\end{proof}

\begin{lemma}\label{lem:quotient-law}
Let $A$ be a vector space, and consider subspaces $B, C, D \subseteq A$ such that $C \subseteq B$. Then, there is an isomorphism
\[
\dfrac{B+D}{C+D}
\cong 
\dfrac{B}{C+ B \cap D} .
\]
\end{lemma}

\begin{proof}
Consider the inclusion $\iota\colon B \hookrightarrow B + D$, which induces a linear map
\[
\begin{tikzcd}[/tikz/column 1/.style={column sep=-0.5em}]
\oviot : & 
\dfrac{B}{C+ B \cap D} \arrow[r] & 
\dfrac{B+D}{C+D}.
\end{tikzcd}
\]
First, notice that $\oviot$ is well-defined since $\iota(C+ B \cap D) \subseteq C + D$. We claim that $\oviot$ is an isomorphism. 
To prove injectivity, let $b \in B$ and assume that $\oviot(\overline{b})$ is trivial; that is, $b \in C + D$. Hence, there exists some vector $c \in C$ such that $b - c \in D$. 
However, $b-c \in B$ by the hypothesis $C\subseteq B$. Altogether, $b = c + (b-c) \in C + B\cap D$ and injectivity follows.
Surjectivity follows since, for any class $\overline{b+d}$ in the codomain of $\oviot$, we have the equality $\overline{b+d}=\overline{b}$ and $\overline{b} = \oviot(\overline{b})$.
\end{proof}

\begin{lemma}\label{lem:remove_inter}
    Let $C \subset B$ and $A$ be vector spaces such that $B \subset A + C$, then we have the isomorphism
\[
    \dfrac{B \cap A}{C \cap A} \cong \dfrac{B}{C}
\]
\end{lemma}
\begin{proof}
    The injectivity is clear since $x \in B \cap A$ and $x \in C$ implies $x \in C \cap A$.
    The surjectivity follows from the fact that every $b \in B$ can be written as $a + c$ with $a \in A$ and $c \in C$,
    then $b - c = a$ is in $B$ and $a \in B \cap A$.
    This implies that $a$ has the same class as $b$.
\end{proof}

\section{The induced partial matching}
\label{sec:induced}

In this section, we define the induced partial matching $\matching{f}$ using a variant of a matrix reduction algorithm.
Throughout this section, we fix two persistence modules, $V$ and $U$, and a basis for each: $\cA = \{\alpha_j\}_{j\in [\barsize{V}]}$, and $\cB = \{\beta_i\}_{i\in [\barsize{U}]}$.
In order to perform the matrix reduction, we need to define a matrix representation of a persistence map $f \colon V \to U$.
\begin{definition}\label{def:assoc-matrix}\citep{distributed}
    The \emph{matrix} $F$ associated with $f$ with respect to $\cA$ and $\cB$ is a $\barsize{U} \times \barsize{V}$ matrix with coefficients in $k$ such that, for each generator $\alpha_j \in \cA$,
    \[
        f\big(\alpha_{ja_j}^1\big) = \sum_{\beta_{ia_j} \in \cB_{a_j}} F(i,j) \beta_{ia_j}^1.
    \]
    where $i \in [\barsize{U}]$, $j \in [\barsize{V}]$ and $F(i,j)$ are the entries of $F$.  
\end{definition}

Note that $F(i,j)=0$ for all $\beta_i \in \cB$ such that $[c_i, d_i]\not\leq [a_j, b_j]$ due to the commutativity of $f$ with the structure maps in (\ref{eq:persistence_map}).
Recall that $\cB_{a_j}$ provides a basis for $U_{a_j}$.
More generally, we define the submatrix of $F$, $F_t$, with columns $j$ such that $t \in [a_j, b_j]$ and rows $i$ with $t \in [c_i, d_i]$.
It is easy to check that $F_t$ is precisely a matrix for the linear map $f_t:V_t \rightarrow U_t$. 
\smallskip

In order to calculate $\matching{f}$, we need to set a specific order on the entries of $F$.
\begin{definition}\label{def:orders}
    We write $[c,d] \leqU [a,b] $ if $d<b$ or whenever $d=b$ and $c\leq a$.
    Analogously, we write $[c,d] \leqV [a,b]$ if $c<a$ or whenever $c=a$ and $d\leq b$.
\end{definition}

The order $\leqU$ was already defined in \citep{distributed} as the \emph{endpoint order}.
Note that $J\leq I$ implies both $J \leqU I$ and $J \leqV I$. 
Throughout the rest of this article, we assume the elements of $\cA$ are ordered such that $j' \leq j$ implies $[a_{j'}, b_{j'}] \leqV [a_j, b_j]$.
Similarly, for $[c_{i'}, d_{i'}],[c_i, d_i] \in \cB$, if $i' \leq i$ then $[c_{i'}, d_{i'}] \leqU [c_i, d_i]$. 
In the following, we assume that the rows and columns of $F$, and the elements of $\cA$ and $\cB$ are ordered in this way.
\smallskip

Before introducing the algorithm, we define the function $\pivot_F$ as
\[
    \pivot_F(i,j)  = 
    \begin{cases}
        1 & \text{if $F(i,j)$ is the lowest non-zero entry of column $j$}\\
        0 & \text{otherwise}.
    \end{cases}
\]
If $\pivot_F(i,j) = 1$, we say that $(i,j)$ is a pivot.
Fixing $\beta_i \sim I$ and $\alpha_j \sim J$, we define $\delta_{ij} = 1$ if $I \cap J \neq \emptyset$ and $\delta_{ij} = 0$ otherwise.
The algorithm $\matching{f}$ reduces $F$ to a new matrix, $R$, using column elimination, see Algorithm~\ref{alg:column_reduction}.
\begin{algorithm}[hbt]
  \caption{Column reduction}\label{alg:column_reduction}
  \begin{algorithmic}[1]
    \STATE \textbf{Input: $F$ (see Definition~\ref{def:assoc-matrix}) }
    \STATE \textbf{Output:} reduced matrix by left-to-right column operations.
	\STATE{$R\gets F$};
    \STATE {$i \gets \barsize{U} + 1 $;}
	\WHILE {$i \neq 0$} 
            \STATE {$i\gets i-1$;}
            \STATE {$A \gets \{ j  \mid \pivot_R(i,j) = 1\}$;}
            \IF{$A \neq \emptyset$}
                \STATE{ $\ell \gets \min A$}
        		\WHILE{$\exists j > \ell$ such that $\pivot_R(i,j) = 1$} 
                    \FOR{$z = 1, \ldots, i$}
                        \STATE{$R(z,j) \gets \delta_{z j}\left( R(z,j)  - \dfrac{R(i,j)}{R(i,\ell)}R(z,\ell)\right)$}
                    \ENDFOR
                \ENDWHILE;
            \ENDIF;
	\ENDWHILE.
    \STATE \textbf{Returns:} $R$
  \end{algorithmic}
\end{algorithm}

After obtaining $R$ with Algorithm~\ref{alg:column_reduction}, we define $\matching{f} \colon\II(n) \times \II(n) \to \mathbb{Z}_{\geq 0}$ as:
\begin{equation}\label{eq:definition}
    \matching{f}(I,J) =
    \#\left\{ 
    (i,j) \text{ pivots of } R \mid [a_j, b_j] = I, [c_i, d_i] = J
    \right\}.
\end{equation}

\begin{theorem}\label{thm:image}
    $\matching{f}$ is a partial matching.
\end{theorem}
\begin{proof}
    We need to prove that $\matching{f}(I,J)$ satisfies Inequalities~(1) and (2) from Definition~\ref{def:partial_matching}.
    We start with (2).
    Note that, by definition,
    \[
    \sum_{I \in \II(V)} \matching{f}(I, J) = \#\left\{ 
    (i,j) \text{ pivots of } R \mid [c_i, d_i] = J
    \right\}.
    \]
    We use now that Algorithm~\ref{alg:column_reduction} assigns at most one pivot to every row $i$.
    To see this, observe that when the iteration arrives at $i$, 
    at most one pivot is left for the row $i$; therefore, at most one column remains with pivot at row $i$.
    This column will be modified again only if it had another pivot in row $i'$ with $i' < i$. However this is impossible since, by definition, there can be at most one pivot per column. 
    Hence,
    \[
    \sum_{I \in \II(V)} \matching{f}(I, J) \leq  \#\left\{i \mid \beta_i \sim J \right\}.
    \]
    Lastly,
    by definition of persistence basis, the number of rows $i$ with $[c_i, d_i] = J$ corresponds to the number of elements in $\cB$ with $\beta_i \sim J$, which is precisely the multiplicity $m_J^U$.
    Then, $ \#\left\{i \mid \beta_i \sim J \right\} = m_J^U$ and the result follows.
    The other inequality can be proven analogously.
\end{proof}
Actually, Algorithm~\ref{alg:column_reduction} calculates the image of $f$ as a byproduct.
This might be intuitive, but it is not straightforward to prove.
The main problem is that column reduction of persistence maps is beset by mathematical difficulties \citep{distributed, matrix}.
Actually, most algorithms for calculating $f V$ are designed in the specific case of persistence modules coming from persistent homology, so they are defined only for chain complexes \citep{cok, efficient_image}.
We will give a precise statement of this fact in Section~\ref{sec:properties}.

\begin{example}\label{exa:calculation_limitation}
    We calculate $\matching{f}$ for $f$ from Example~\ref{exa:limitation}. 
    We use $\alpha_1 \sim [2,2]$ and $\alpha_2 \sim [2,3]$ as explained in Example~\ref{ex:bases_alpha}.
    Fixing the canonical basis for $U$ with $\beta_1 \sim [1,2]$ and $\beta_2 \sim [1,2]$, we get the following $F$,
    \[
    F= \;\;
            \bordermatrix{ & \alpha_1 & \alpha_2 \cr
              \beta_1 & 0 & 0 \cr
              \beta_2 & 1 & 0  }
    \]
    which is already reduced.
    Then, the only non-null value of the induced partial matching is $\matching{f}([2,2], [1,2]) = 1$. \exend
\end{example}

\begin{example}\label{exa:indecomposable_matching}
    Let us consider Example~\ref{ex:complex_indecomposable}. In this case, $V$ has a persistence basis given by two generators $\alpha_1\sim [1,4]$ and $\alpha_2\sim [2,3]$ while $U$ has a persistence basis formed by $\beta_1\sim [2,2]$ and $\beta_2 \sim [1,3]$. 
    The matrix $F$ associated with $f\colon V\rightarrow U$ and its reduced matrix are 
    \[
    F= \;\;
            \bordermatrix{ & \alpha_1 & \alpha_2 \cr
              \beta_1 & 0 & 1 \cr
              \beta_2 & 1 & 1  }\, 
    \longrightarrow 
    R = \;\;
            \bordermatrix{ & \alpha_1 & \alpha_2 \cr
              \beta_1 & 0 & 1 \cr
              \beta_2 & 1 & 0  } .
    \]
    Reading the pivots from $R$, we get
    \begin{align*}
    \matching{f}([1,4], [1,3]) &= \matching{f}([2,3], [2,2])= 1, \text{ and}
    \\
    \matching{f}([1,4], [2,2]) &= \matching{f}([2,3], [1,3])= 0.
    \tag*{\exend}
    \end{align*}
\end{example}

Before finishing this section, we analyze the complexity of Algorithm~\ref{alg:column_reduction}. 
In particular, we show that it is cubic on $M = \max(\#\PD(V), \# \PD(U))$.
To start, the outer \texttt{while} loop is executed $\leq M$ times. 
Next, the complexity of computing $A$ is bounded by $M^2$; assuming that pivots of $M$-length columns are computed ad-hoc.
Also, computing $\min A$ takes $M$ operations at most.
Finally, the inner \texttt{while} loop is executed at most $M$ times, and each time it performs operations on columns of length $\leq M$.
Altogether, we obtain a complexity of $O\left(M( M^2+M+M^2)\right) \sim O\left( M^3 \right)$.

\section{Properties of the induced matching}\label{sec:properties}

The definition of $\matching{f}$ may suggest that it depends on the choice of bases.
To prove that this is not the case, we provide an alternative definition of $\matching{f}$ in terms of persistence modules, images and kernels.
As a direct consequence of this definition, we obtain the additivity of $\matching{f}$ with respect to direct sums.
We first define some operators to this end.
\begin{definition}\label{def:V-tilde-hat}
 For a fixed interval $[a,b]$, and a persistence module $V$ with structure maps $\rho$, we define the following operators,
\begin{align*}
\tV_{[a,b]t} &\coloneqq \begin{cases*}
                V_t & if  $t \in [1,a-1]$, \\
                \im \rhov{a}{t} \cap \ker \rhov{t}{(b+1)} + \im \rhov{(a-1)}{t} & if  $t \in [a,b]$,  \\
                 \im \rhov{(a-1)}{t} & if  $t \in [b + 1,n]$.
            \end{cases*}
\\
\hV_{[a,b]t} &\coloneqq \begin{cases*}
                \ker \rhov{t}{(b+1)} & if  $t \in [1,a-1]$,  \\
                \im \rhov{a}{t} \cap \ker \rhov{t}{(b+1)} + \ker \rhov{t}{b} & if  $t \in [a,b]$,  \\
                 0 & otherwise.
            \end{cases*}
\end{align*}
\end{definition}

From Lemma~\ref{lem:commute_intersection}, we have that $\rhov{s}{t} \tV_{[a,b]s} \subset \tV_{[a,b]t}$ and $\rhov{s}{t} \hV_{[a,b]s} \subset \hV_{[a,b]t}$.
Then, they are persistence modules whose structure map is the restriction of $\rho$.
As the following lemma shows, the nature of these operators become much simpler when expressed in terms of persistence bases.
\begin{lemma}\label{lem:Abasis}
    Given a persistence basis $\cA$ of $V$, and an interval $[a,b]$, we define
\begin{align*}
    \tcA_{[a,b]}&\coloneqq\big\{ 
        \alpha_i \in \cA \mid (a_i < a) \text{ or } (a_i=a \text{ and } b_i\leq b)
    \big\},
    \\
    \hcA_{[a,b]}&\coloneqq \big\{ 
        \alpha_i \in \cA \mid (b_i < b) \text{ or } (b_i=b \text{ and } a_i\leq a)
    \big\}, \mbox{ and }
    \\
    \cA_{[a,b]}&\coloneqq \big\{ 
        \alpha_i \in \cA \mid (a_i=a \text{ and } b_i=b)
    \big\}.
\end{align*}
Then,
\begin{itemize}
    \item $\tcA_{[a,b]}$ is a persistence basis of \,$\tV_{[a,b]}$, 
    \item $\hcA_{[a,b]}$ is a persistence basis of \,$\hV_{[a,b]}$, and
    \item $\cA_{[a,b]}$ is a persistence basis of $\tV_{[a,b]}\big/\tV_{[a,b-1]}$ and $\hV_{[a,b]}\big/\hV_{[a-1,b]}$.
\end{itemize}
\end{lemma}
\begin{proof}
    In \citep[Lemma~3.1]{InducedMatchings2022}, it was proven that $\im \rhov{a}{t}$ is generated by $\cI_{at}$, where
    \[
        \cI_{a} \coloneqq
        \left\{\alpha_i \in \cA \mid \alpha_i \sim [a_i, b_i] \text{ such that } a_i \leq a \right\}
    \]
    and analogously, 
    $\ker \rhov{t}{(b+1)}= \langle \cK_{bt}\rangle$, where 
    \[
        \cK_{b} \coloneqq
        \left\{ \alpha_i \in \cA \mid \alpha_i \sim [a_i, b_i] \text{ such that } b_i \leq b \right\}.
    \]
    Hence, using these results and the definition of $\tV_{[a,b]}$, we get that for any $t\in [n]$, $\langle \tcA_{[a,b]t} \rangle$ is isomorphic to $\langle \cI_{at} \cap \cK_{bt} \cup \cI_{(a-1)t} \rangle$, which coincides with $\tcA_{[a,b]t}$.
    The other cases are analogous.
\end{proof}

\begin{example}
    Consider $V$ from Example~\ref{exa:diagram}, and recall that it is isomorphic to $k_{[2,2]}\oplus k_{[2,3]}$. 
    Hence, $\tV_{[2,3]}$ is isomorphic to $V$, $\tV_{[2,2]}$ is isomorphic to $k_{[2,2]}$, and $\tV_{[2,3]}/\tV_{[2,2]}$ is isomorphic to $k_{[2,3]}$. \exend
\end{example}

From Lemma~\ref{lem:Abasis}, we get an alternative way of calculating the decomposition of $V$. 
For any $t\in[a,b]$,
\[
    \dim \tV_{[a,b]t}\big/\tV_{[a,b-1]t} = \# \cA_{[a,b]} = m_{[a,b]}^V.
\]
In addition, since the generators commute with the structure maps, 
we also have that  
\[
\tV_{[a,b]} \big/ \tV_{[a,b-1]} \cong \hV_{[a,b]}\big/\hV_{[a-1,b]} \cong  \bigoplus_{1}^{m^V_I}k_{[a,b]}.
\]
\begin{definition}
    For every persistence module $V$ and interval $[a,b]$, we define the persistence module
    \[
        V_{[a,b]} \coloneqq \tV_{[a,b]} \big/ \tV_{[a,b-1]}.
    \]
\end{definition}
Note that $V_{[a,b]}$ is also isomorphic to $\hV_{[a,b]}\big/\hV_{[a-1,b]}$ and has $\cA_{[a,b]}$ as a persistence basis.
One might wish to combine $V_I$ and $U_J$ to find the relation with $\matching{f}(I,J)$.
However, since they are quotients of different persistence modules, $V_I \cap U_J$ is not well defined. 
Instead, we provide the following alternative definition.

\begin{definition}\label{def:Zpm}
    Given two intervals $I=[a,b]$ and $J=[c,d]$, and a persistence map $f: V \to U$, we define the following vector spaces for $t\in[n]$,
    \[
        \subU_{IJt}(f) \coloneqq
         \dfrac{f \tV_{[a,b]t} \cap \hU_{[c,d]t}}{f\tV_{[a,b-1]t} \cap \hU_{[c,d]t} + f \tV_{[a,b]t}\cap \hU_{[c-1,d]t}} .
    \]
    which form a persistence module, $\subU_{IJ}(f)$,  with structure maps the restriction of $\varphi$.
\end{definition}
 Note that we deduce directly that $\subU_{IJ}(f)$ is a persistence module, as it consists of the intersection, sum and quotient of submodules of $U$.
 Next, we present a matrix $R^{IJ}$ such that $\subU_{IJt}=\langle R^{IJ}_t\rangle$ for all $t \in I\cap J$; where by $\langle R^{IJ}_t\rangle$ we denote the vector space generated by the columns from $R^{IJ}_t$.
\begin{definition}
    Let $f : V \to U$ be a persistence map, $\cA$ and $\cB$ persistence bases for $V$ and $U$, and $R$ the reduced matrix given by Algorithm~\ref{alg:column_reduction}.
    Denote the rows $i$ with $[c_i, d_i] = J$ as $\cB_J$ and the columns $j$ with $[a_j, b_j] = I$ as $\cA_I$.
    Then, $R^{IJ}$ is the submatrix of $R$ with rows in $\cB_J$ and columns from $\cA_I$ containing pivots in $\cB_J$.
\end{definition}

\begin{theorem} \label{the:ZandM}
We have that $\langle R^{IJ}_t \rangle \cong \subU_{IJt}(f)$ and 
    \[
        \matching{f}(I,J) = \rank \langle R^{IJ}_t \rangle = \dim \subU_{IJt}(f)
    \]
    for any $t \in I \cap J$.
\end{theorem}
\begin{proof}
    The first equality is straightforward by definition of $\matching{f}(I,J)$ and $R^{IJ}$.
    For the second one, define the matrix 
    $Y^{IJ}$ as the submatrix of $R$ formed by the columns from $\tcA_{I}$ with pivots in $\cB_{J}$.
    Let $[c,d]$ be the interval $J$.
    Notice that $\langle Y^{IJ}_t \rangle$ is isomorphic to the vector space generated by the columns from $\tcA_{I}$ with pivots in $\hcB_{[c,d]}$, modulo its subspace generated by columns of $\tcA_{I}$ with pivots in $\hcB_{[c-1, d]}$.
    More precisely,
    \[
        \langle Y_{t}^{IJ} \rangle \simeq \dfrac{f(\langle \tcA_{It} \rangle) \cap \langle \hcB_{[c,d]t} \rangle}{f(\langle \tcA_{It} \rangle) \cap \langle \hcB_{[c-1,d]t} \rangle}.
    \]
    Define
    \begin{align*}
    B &= f(\langle \tcA_{It} \rangle) \cap \langle \hcB_{[c,d]t} \rangle, 
    \\
    C &= f(\langle \tcA_{It} \rangle) \cap \langle \hcB_{[c-1,d]t} \rangle, \text{ and}
    \\
    D &= \langle \hcB_{[c-1,d]t} \rangle;
    \end{align*}
    and notice that $C = B \cap D$.
    Hence, we can apply Lemma~\ref{lem:quotient-law}, to get
    \[
        \langle Y_{t}^{IJ} \rangle \simeq \dfrac{f(\langle \tcA_{It} \rangle) \cap \langle \hcB_{[c,d]t} \rangle +  \langle \hcB_{[c-1,d]t} \rangle}{\langle \hcB_{[c-1,d]t} \rangle}.
    \]
    We now write $I$ as $[a,b]$.
    It is clear from the definition that 
    \[
        \langle R^{[a,b]J}_t \rangle \cong \langle Y^{[a,b]J}_t \setminus Y^{[a,b-1]J}_t \rangle \cong \langle Y^{[a,b]J}_t \rangle / \langle  Y^{[a,b-1]J}_t \rangle.
    \]
    Hence,
    \[
        \langle R^{[a,b]J}_t \rangle 
        \cong 
        \dfrac{
        f(\langle \tcA_{[a,b]t} \rangle) \cap \langle \hcB_{[c,d]t} \rangle +  \langle \hcB_{[c-1,d]t} \rangle 
        }{
        f(\langle \tcA_{[a,b-1]t} \rangle) \cap \langle \hcB_{[c,d]t} \rangle +  \langle \hcB_{[c-1,d]t} \rangle
        }
    \]
    Using Lemma~\ref{lem:quotient-law} with 
    \begin{align*}
        B &= f(\langle \tcA_{[a,b]t} \rangle) \cap \langle \hcB_{[c,d]t} \rangle,
        \\
        C &= f(\langle \tcA_{[a,b-1]t} \rangle) \cap \langle \hcB_{[c,d]t} \rangle, \text{ and} 
        \\
        D &= \hcB_{[c-1,d]t};
    \end{align*}
    we get that $\langle R^{[a,b]J}_t \rangle$ is isomorphic to
    \begin{equation}\label{for:isoZR}
    \dfrac{f(\langle \tcA_{[a,b]t} \rangle) \cap \langle \hcB_{[c,d]t} \rangle}{f(\langle \tcA_{[a,b-1]t} \rangle) \cap \langle \hcB_{[c,d]t} \rangle + f(\langle \tcA_{[a,b]t} \rangle) \cap \langle \hcB_{[c-1,d]t} \rangle}
    \end{equation}
    which is itself isomorphic to the vector space $\subU_{IJt}(f)$ by Definition~\ref{def:Zpm}.
\end{proof}
There are two important immediate consequences of this result:
$\matching{f}$ is well-defined and does not depend on the bases chosen for Algorithm~\ref{alg:column_reduction} and,
$\matching{f}$ is additive with respect to direct sums.

\begin{corollary}
    $\matching{f} \colon\II(n) \times \II(n) \to \mathbb{Z}_{\geq 0}$ is an invariant of $f$.
\end{corollary}
\begin{proof}
    By definition, if $f \cong g$ then $\subU_{IJ}(f) \cong \subU_{IJ}(g)$, and
    by Theorem~\ref{the:ZandM}, 
    \[
        \matching{f}(I,J) = \dim \subU_{IJt}(f) = \dim \subU_{IJt}(g) = \matching{g}(I,J). \qedhere
    \]
\end{proof}

\begin{corollary}\label{cor:additive}
    Given two persistence maps,
    $f\colon V \rightarrow U$ and 
    $g \colon V' \rightarrow U'$,
    and intervals $I,J \in \II(n)$, we have that
    \[
        \matching{f \oplus g} (I, J) = 
        \matching{f}(I,J) + \matching{g}(I,J).
    \]
\end{corollary}
\begin{proof}
    The result follows since direct sums commute with quotients, finite intersections and sums of vector spaces, and $\dim (V_t \oplus U_t) = \dim V_t + \dim U_t$.
    Hence,
    \[ \matching{f \oplus g} (I, J) = \dim \left( \subU_{IJt}(f \oplus g) \right) = \dim \subU_{IJt}(f) + \dim \subU_{IJt}(g). \qedhere \]
\end{proof}

The previous corollary indicates that $\matching{f}$ depends on the decomposition of $f$. 
For instance, in Example~\ref{exa:limitation}, the only non-null value of $\matching{f}$ is 
$\matching{f}([2,2], [1,2]) = 1$.
\smallskip

Theorem~\ref{the:ZandM} also characterizes the  structure of $\subU_{IJ}$.
As the following result shows, it is formed by the interval modules obtained after intersecting the intervals matched by $\matching{f}$.
\begin{corollary}\label{cor:decomposition_Z}
    Given any pair $s\leq t$ from $I\cap J$, we have that $\phiu{s}{t}$ induces an isomorphism $\subU_{IJs}(f)\cong \subU_{IJt}(f)$.
    Moreover,
    \[
        \subU_{IJ}(f) \cong \bigoplus_{1}^{\matching{f}(I,J)} k_{I \cap J} .
    \]
\end{corollary}
\begin{proof}
    Note that $\rhov{s}{t}$ induces an isomorphism between $\langle    R^{IJ}_s\rangle$ and $\langle R^{IJ}_t \rangle$ for $s,t \in I \cap J$.
    In addition, this isomorphism commutes with the one of Expression~(\ref{for:isoZR}).
    Lastly, using Definition~\ref{def:V-tilde-hat}, one can check that $\subU_{IJ}(f)\cong 0$ for all $t \notin I \cap J$.
    For instance, if $I=[a,b]$ and $t<a$, then $\tV_{[a,b]t} = \tV_{[a,b-1]t} = V_t$ and so (using $\hU_{[c-1,d]t} \subset \hU_{[c,d]t}$) it follows that $\subU_{IJ}(f)=(fV_t\cap \hU_{[c,d]t}) / (fV_t\cap \hU_{[c,d]t}) \cong 0$.
\end{proof}
Since the definition of $\subU_{IJt}$ can be extended to continuous persistence modules, it opens the door to a generalization to continuous persistence modules of $\matching{f}$ itself, similarly to the constructions in \citep{InducedMatchings2022}.

\subsection{The induced projection}\label{sec:projection}
There is also an interesting interpretation of $\matching{f}$ as a projection in the space of persistence maps.
First, we write $k_{IJ}\colon k_I \to k_J$ for the persistence map that is zero unless $J\leq I$; in which case it is the identity between $k_{It}$ and $k_{Jt}$ for any $t \in I \cap J$.
The induced partial matching behaves as expected with respect to these maps.
\begin{proposition}\label{pro:interval-case}
    Let $J\leq I$ and
    $k_{IJ} \colon k_I \to k_J$. The only non-null value of $\matching{k_{IJ}}$ is $\matching{k_{IJ}}(I,J) = 1$.
\end{proposition}
\begin{proof}
    The result follows directly from Formula~(\ref{eq:definition}), since it is a matrix with only one entry.
    The column corresponds to the interval $J$, the row to $I$, and the entry is non-null, so it is a pivot.
\end{proof}

It is straightforward to check that non-null maps between interval modules are indecomposable \citep{ladder}.
We say that a persistence map is \emph{interval decomposable} if $f \cong f_1 \oplus \ldots \oplus f_n$ and every $f_i$ is a non-null map between interval modules.
Hence, we can define an operator sending a persistence map, $f \colon V \to U$, to interval decomposable maps using $\matching{f}$:
\[
    {\it P}(f) = \bigoplus_{\substack{I \in \II(V) }}
    \bigoplus_{\substack{J \in \II(U) }} \left( \bigoplus_{1}^{\matching{f}(I,J)} k_{IJ} \right).
\]
\begin{corollary}
    The operator ${ \it P }$ is an additive projection of persistence maps to interval decomposable maps.
\end{corollary}
\begin{proof}
    The additivity is a direct consequence of Corollary~\ref{cor:additive}.
    We also have that ${\it P}^2 = {\it P}$ since ${\it P} (k_{IJ}) = k_{IJ}$ by Proposition~\ref{pro:interval-case}.
\end{proof}

\subsection{Relation with $\chi_f$}

As mentioned in the introduction, the induced partial matching $\chi_f$ has been extensively used,
so it is natural to ask whether it is related to $\matching{f}$.
From Example~\ref{exa:limitation}, we can deduce that $\chi_{f}$ is not additive and is hence different from $\matching{f}$. 
However, there is more to be said.
To avoid a significant detour from the main goals of this article, 
we present a complete comparison between the two matchings in \ref{app:BL}, where the formal definition of $\chi_f$ is recalled.
As a summary, we highlight the following results:
\begin{itemize}
    \item Do $\chi_f$ and $\matching{f}$ coincide for indecomposable persistence maps? In other words, is $\matching{f}$ an ``additivization'' of $\chi_f$? The answer is no, as shown by Example~\ref{ex:m_isnot_chi}.
    \item Given $f \colon V \to U$, $\chi_f$ is completely defined by $f V$.
    However, we prove in Section~\ref{sec:image} that $f V$ can be obtained from $\matching{f}$, and moreover, $\matching{f}$ is a richer invariant than $f V$ (see Example~\ref{exa:same_image}).
    In particular, $\matching{f}$ is richer than $\chi_f$.
    \item There exists a class of persistence maps for which their internal structure is greatly simplified, making $\matching{f}$, $\chi_f$, and $f V$ coincide.
    We call them \emph{polybar} maps and introduce them in \ref{app:BL}.
\end{itemize}

\section{Obtaining the image module}\label{sec:image}
In this section, we prove that $\matching{f}$ contains strictly more information than $f V$ and that, as a consequence, Algorithm~\ref{alg:column_reduction} calculates the image module.
More concretely, we prove the following statement.

\begin{theorem}\label{the:image}
Given $f\colon  V \rightarrow U$, we have the isomorphism
    \[ 
    f V 
    \cong 
    \bigoplus_{\substack{I \in \II(V) }}
    \bigoplus_{\substack{J \in \II(U) }}
    \subU_{IJ}(f).
    \]  
\end{theorem}
We prove this result in Section~\ref{sec:proof}, after introducing the main technical tool in Section~\ref{sec:sections}.
Before doing so, note that we obtain a direct way of calculating the decomposition of $f V$ from $\matching{f}$.
\begin{corollary}\label{cor:mult-intervals-im}
    Given $L \in \II(f V)$,
    we have that
    \[ 
        m_L^{fV} = \sum_{\substack{I \cap J = L}} \matching{f}(I,J). 
    \]
    In particular, $\PD(f V)$ can be calculated from $\matching{f}(I,J)$.
\end{corollary}
\begin{proof}
    By Theorem~\ref{the:image}, each summand $\subU_{IJ}$ contributes $\matching{f}(I,J)$ copies of $k_{I\cap J}$ in the decomposition of $f V$.
    Thus, the number of copies of $k_{L}$ is determined by the sum of $\matching{f}(I,J)$ such that $I \cap J = L$.
\end{proof}
The following example illustrates the fact that $\matching{f}$ determines $fV$.
\begin{example}
    Recall $f$ from Example~\ref{ex:complex_indecomposable}.
    The induced partial matching was calculated in Example~\ref{exa:indecomposable_matching}:
    \begin{align*}
        \matching{f}([1,4], [1,3]) &= \matching{f}([2,3], [2,2])= 1, \text{ and}
        \\
        \matching{f}([1,4], [2,2]) &= \matching{f}([2,3], [1,3])= 0. 
    \end{align*}
    Then, the intersections of intervals with non-null $\matching{f}$ are $[1,4] \cap [1,3] = [1,3]$ and $[2,3] \cap [2,2] = [2,2]$.
    Hence, the persistence diagram of $f V$ is
    $\{([1,3], 1), ([2,2],1) \}$.
    \exend
\end{example}
Moreover, as the following example shows, $\matching{f}$ is a finer invariant than the image module. 
\begin{example}\label{exa:same_image}
The following persistence module, $g$, has the same image, domain and codomain modules as $f$ from Example~\ref{exa:limitation}, but a different induced partial matching (cf. Example~\ref{exa:calculation_limitation}):
 \begin{equation*}	
    \begin{tikzcd}[cramped]
        k \arrow[r, "\Id"] & k \arrow[r] & 0 \\
        0 \arrow[r]\arrow[u] & 0 \arrow[r]\arrow[u] & 0\arrow[u]
    \end{tikzcd}
        \,\oplus\,
    \begin{tikzcd}[cramped]
        0 \arrow[r] & 0 \arrow[r] & 0 \\
        0 \arrow[r]\arrow[u] &  k \arrow[r]\arrow[u, "\Id"] & 0 \arrow[u]
    \end{tikzcd}
        \,\oplus\,
    \begin{tikzcd}[cramped]
        k \arrow[r, "\Id"] & k \arrow[r] & 0 & \\
        0 \arrow[r]\arrow[u] & k \arrow[r, "\Id"]\arrow[u] & k. \arrow[u] 
    \end{tikzcd} 
\end{equation*}
By Corollary~\ref{cor:additive}, the only non-null relation is $\matching{g}([2,3], [1,2]) = 1$.
\exend
\end{example}

In order to prove Theorem~\ref{the:image}, we need to introduce sections, a tool for decomposing persistence modules.
They are explained in the next subsection.

\subsection{Sections of a vector space}\label{sec:sections}

Sections were first introduced in TDA in~\citep{CrawleyBoevey2015} to prove Theorem~\ref{the:decomposition}.
A \emph{section} of a vector space 
$A$ 
is a pair of vector spaces 
$(F^-,F^+)$
such that $F^- \subseteq F^+ \subseteq
A$. 
We say that a set $\{ (F^-_\lambda,F^+_\lambda)\colon\lambda\in \Lambda \}$  of sections of $A$ with index set $\lambda$ is \emph{disjoint} if, for all $\lambda\neq \mu$,
either $F^+_\lambda \subseteq F^-_\mu$ or $F^+_\mu \subseteq F^-_\lambda$.
Also, such a set \emph{covers $A$} provided that for any subspace 
$B \subsetneq A$
there is some $\lambda \in \Lambda$ such that
$    B + F^-_\lambda \neq 	B + F^+_\lambda $.

\begin{theorem}\citep[Theorem~6.1]{CrawleyBoevey2015}\label{the:subsumcond}
	Let $A$ be a vector space and let $\{(F^-_\lambda,F^+_\lambda)\colon\lambda\in \Lambda \}$ be a disjoint set of
	sections that covers $A$.
    Then,
    \[
    \bigoplus_{\lambda\in\Lambda} 
    	    \big( F^+_\lambda\, \big/ \,F^-_\lambda\big)\cong A.
    \]
\end{theorem}

In general, finding a set of sections that covers $A$ is not an easy task. 
Alternatively, one can combine simpler covers to obtain a new one.
A set of sections $\{(F^-_\lambda,F^+_\lambda)\colon\lambda\in \Lambda\}$ \emph{strongly covers} a vector space $A$ provided that, for all subspaces $B, C\subseteq A$
with $C\not\subseteq B$, there is some $\lambda \in \Lambda$ with
\[
B+(F^-_{\lambda}\cap C)\neq 
B+(F^+_{\lambda}\cap C).
\]
Strongly covering sets of sections of $A$ can be combined to produce a new one, as shown by the following lemma.
\begin{lemma}\citep[lemma~6.2]{CrawleyBoevey2015}\label{lem:secmix}
    If $\{(F^-_\lambda, F^+_\lambda)\colon\lambda \in \Lambda \}$ is a disjoint set of sections that  covers a vector space $A$, and $\{(G^-_\sigma, G^+_\sigma):\sigma \in \Sigma\}$ is a disjoint set of sections that strongly covers $A$, then the set
    $ \{( G^-_{\sigma}\cap F^+_{\lambda}+F^-_{\lambda},\;  G^+_\sigma \cap F^+_\lambda+F^-_{\lambda} )\colon(\lambda, \sigma) \in \Lambda \times \Sigma \},
    $
    is a disjoint set of sections that covers $A$.
\end{lemma}
    Lastly, we include a technical lemma that will be needed later.
    \begin{lemma}\label{lem:3cond}
        Consider a vector space $A$, a finite totally ordered set $\{\theta_i\}_{i \in [m]}$, and a set of sections of $A$, $\left\{F_{\theta_i}^-, F_{\theta_i}^+\right\}_{i \in [m]}$.
        If the following conditions hold
            \begin{enumerate}
                \item $F_{\theta_i}^+ = F_{\theta_{i+1}}^-$ for all $1 \leq i < m$,
                \item $F^-_{\theta_1} = 0$, and 
                \item $F^+_{\theta_m} = A$;
            \end{enumerate}
        then $\left\{F_{\theta_i}^-, F_{\theta_i}^+\right\}_{i \in [m]}$ is a disjoint set of sections that strongly covers $A$.
\end{lemma}
\begin{proof}
We start proving the disjoint property.
Let $i < j$,
by the first condition,
there exists a sequence of spaces
$
    F^+_{\theta_i} = F^-_{\theta_{i + 1}} \subseteq F^+_{\theta_{i + 1}} = \ldots = F^-_{\theta_{j}}
$
and $F^+_{\theta_i}  \subseteq F^-_{\theta_j}$.
\smallskip

Now, let us prove that
$\{F_{\theta_i}^-, F_{\theta_i}^+\}_{i \in [m]}$ strongly covers $A$. For any $C \nsubseteq B$ with $C,B \subseteq A$, define the sets
\[ 
\Theta_{B,C}^- = \{ \theta_i \, \mid \, F^-_{\theta_i}\cap C \subseteq B \},
\quad 
\Theta_{B,C}^+ = \{ \theta_i \, \mid \, F^+_{\theta_i}\cap C \nsubseteq B\},
\]
and note that due to conditions $2$ and $3$ from the lemma statement, 
neither of these sets is empty.
In addition, if $\theta_l = \max \{\Theta_{B,C}^- \}$, then $\theta_{l + 1} \notin  \Theta_{B,C}^-$. 
Since $F^+_{\theta_l} = F^-_{\theta_{l+1}}$, $\theta_l$ must be in   $\Theta_{B,C}^+$ and
\[
    B + F_{\theta_l}^- \cap C 
    =
    B
    \neq
    B + F_{\theta_l}^+\cap C.
\]
Note that this result also holds when $\theta_l = \theta_m$.
Hence, $\left\{F_{\theta_i}^-, F_{\theta_i}^+\right\}_{i \in [m]}$ strongly covers $A$.
\end{proof}

\begin{remark}\label{rem:subindices}
    Notice that if $\{F_\lambda^-, F_\lambda^+\}_{\lambda \in \Lambda}$ is a disjoint set of sections of $A$ with $\Omega \subseteq \Lambda$,
    and $\{F_\omega^-, F_\omega^+\}_{\omega \in \Omega}$ (strongly)
    covers $A$, then so
    does $\{F_\lambda^-, F_\lambda^+\}_{\lambda \in \Lambda}$.
\end{remark}

\subsection{The proof of Theorem~\ref{the:image}} \label{sec:proof}
We begin the proof by defining sections for $fV$ and $U$.
\begin{align*}
    f \tV_{[a,b]t}^+ &\coloneqq f (\tV_{[a,b]t}), & \hU_{[c,d]t}^+ &\coloneqq \hU_{[c,d]t}, \\
    f \tV_{[a,b]t}^- &\coloneqq f (\tV_{[a,b-1]t}), & \hU_{[c,d]t}^- &\coloneqq \hU_{[c-1,d]t}.
\end{align*}
We now prove that these sets of sections strongly cover $f V$ and $U$ respectively.
Notice that here is where the orders introduced in Definition~\ref{def:orders} reveal their importance.
\begin{proposition}\label{prop:strong_hV}
    For each $t \in [n]$, 
    we have that 
    \begin{itemize}
        \item[\small $\bullet$]  $\big\{ f \tV_{It}^-, f \tV_{It}^+ \big\}_{I \in \II(n)}$ strongly covers $fV_t$, and
        \item[\small $\bullet$]  $\big\{\hU_{Jt}^-, \hU_{Jt}^+\big\}_{J \in \II(n)}$ strongly covers $U_t$.
    \end{itemize}
\end{proposition}
\begin{proof}
    For a fixed $t$,
    we define the set 
    $\II(n,t) = \{ I \in \II(n) \mid t \in I
    \}$
    and claim that the subset $\big\{ f \tV_{It}^-, f \tV_{It}^+ \big\}_{I \in \II(n,t)}$ strongly covers $f V_t$.
    Then, the result follows by the remark after Lemma~\ref{lem:3cond}.
    We proceed using Lemma~\ref{lem:3cond}.
    First, we use $\leqV$ to fix a total order on the intervals $I \in \II(n,t)$,
    \begin{align*}
    &[1,t] \leqV [1,t+1] \leqV \ldots \leqV [1, n] \leqV [2, t] \leqV [2, t + 1]
    \leqV \ldots \\
    & \ldots  \leqV  [t , n - 1] \leqV [t, n],
    \end{align*}
    Then, since we are assuming that $V_0 = V_{n+1} = 0$, we have that
    \begin{align*}
        &f \tV_{[1,t]t}^- = f (\tV_{[1,t-1]t}) = f \left( \im \rhov{1}{t} \cap \ker \rhov{t}{t} + \im \rhov{0}{t} \right) = 0,\\
        &f \tV_{[t,n]t}^+\supseteq f \left( \im \rhov{t}{t} \cap \ker \rhov{t}{n+1} +  \im \rhov{(t-1)}{t} \right) = f V_t.
    \end{align*}
    These are the two last requirements of Lemma~\ref{lem:3cond}.
    For the first requirement, consider a pair of consecutive intervals in the 
    order $\leqV$; there are two cases:
    \begin{enumerate}[label=(\roman*)]
        \item\label{opt:equal_endpoint} $[a, b-1] \leqV [a,b]$ for some $a,b\in [n]$ such that $t \in 
        [a,b-1]$, or 
        \item\label{opt:change_endpoint} $[a-1,n] \leqV [a,t]$ for some $a\in [n]$ with $a \leq t$.
    \end{enumerate}
    Case \ref{opt:equal_endpoint} follows directly from the definition, since $f \tV^+_{[a, b-1]t} = f \tV^-_{[a, b]t}$.
    For case~\ref{opt:change_endpoint},
    we first notice that $f \tV^+_{[a-1,n]t}$ equals
    \begin{align*} 
     f\left( \im \rhov{(a-1)}{t} \cap V_t + \im \rhov{(a-2)}{t} \right)
    = f\left( \im \rhov{(a-1)}{t} \right),
    \end{align*}
    and that $ \im \rhov{(a-1)}{t} = \tV_{[a,t-1]t} = \tV^-_{[a,t]t}$.
    Altogether, we get $f \tV^+_{[a-1,n]t} = f \tV^-_{[a,t]t}$.
    Hence, we can apply Lemma~\ref{lem:3cond} to obtain that $\big\{\big(f \tV_{It}^-, f \tV_{It}^+\big) \big\}_{I \in \II(n,t)}$ strongly covers $f V_t$. 
    \smallskip
    
    A similar reasoning can be applied to $\big\{ \hU_{Jt}^-, \hU_{Jt}^+ \big\}_{J \in \II(n,t)}$, but using the ordering $\leqU$: 
    \begin{align*}
    &[1,t] \leqU [2,t] \leqU \ldots \leqU [t, t] \leqU [1, t + 1] \leqU [2, t + 1]
    \leqU \ldots \\
    &\ldots \leqU [t - 1, n] \leqU [t, n]. \qedhere
    \end{align*}
\end{proof}
Finally, we have all the ingredients to recover $f V$ from $\subU_{IJt}(f)$.
\begin{proof}[Proof of Theorem~\ref{the:image}]
By Proposition~\ref{prop:strong_hV}, the set of sections $\{\hU_{Jt}^-, \hU_{Jt}^+\}_{J \in \II(n)}$ strongly covers $U_{t}$. Then, using the definition of strongly covering sets of sections, we can take any $B \subsetneq f V_t$ and $C=f V_t$, so that there exists $J\in \II(n)$ such that $B + \hU_{Jt}^- \cap fV_t \neq B + \hU_{Jt}^+ \cap fV_t$. In particular, we deduce that $\big\{ \hU_{Jt}^- \cap fV_t, \hU_{Jt}^+ \cap fV_t \big\}_{J \in \II(n)}$ covers $fV_t$.
In addition,
using Proposition~\ref{prop:strong_hV} and Lemma~\ref{lem:secmix}, 
$W_{IJt}^\pm = f\tV^\pm_{It} \cap \hU^+_{Jt} + \hU^-_{Jt} \cap fV_t$ form a disjoint set of sections that covers $fV_t$.
By Theorem~\ref{the:subsumcond},
    \[
        \bigoplus_{\substack{I \in \II(n) }} \bigoplus_{\substack{J \in \II(n)}} W^+_{IJt} \big/ W^-_{IJt}
        \cong f V_t.
    \]
    Moreover, since $f\tV^\pm_{It} \subset fV_t$, we can use Lemma~\ref{lem:commute_intersection} to obtain
    \[
        W^\pm_{IJt} = \left(f\tV^\pm_{It} \cap \hU^+_{Jt} + \hU^-_{Jt} \right) \cap fV_t.
    \]
    Denoting $A = f V_t$, $ B = f\tV^+_{It} \cap \hU^+_{Jt} + \hU^-_{Jt}$ and $C = f\tV^-_{It} \cap \hU^+_{Jt} + \hU^-_{Jt}$,
    it can be easily seen that $B \subset A + C$.
    Then, we can use Lemma~\ref{lem:remove_inter} to obtain
    \[
        \dfrac{W^+_{IJt}}{W^-_{IJt}}=
       \dfrac{\left(f\tV^+_{It} \cap \hU^+_{Jt} + \hU^-_{Jt} \right) \cap fV_t}{\left(f\tV^-_{It} \cap \hU^+_{Jt} + \hU^-_{Jt} \right) \cap fV_t} 
       \cong 
       \dfrac{f\tV^+_{It} \cap \hU^+_{Jt} + \hU^-_{Jt} }{f\tV^-_{It} \cap \hU^+_{Jt} + \hU^-_{Jt} },
    \]
    Using Lemma~\ref{lem:quotient-law} with $B = f\tV^+_{It} \cap \hU^+_{Jt}$, $C = f\tV^-_{It} \cap \hU^+_{Jt}$ and $D = \hU^-_{Jt}$, we get
    \[
        \dfrac{W^+_{IJt}}{W^-_{IJt}}
        \cong
        \dfrac{f\tV^+_{It} \cap \hU^+_{Jt}}{f\tV^-_{It} \cap \hU^+_{Jt} + f\tV^+_{It} \cap \hU^-_{Jt}}
    \]
    which is precisely the definition of $\subU_{IJt}(f)$.
    Putting all together
    \[
        f V_t
        \cong
        \bigoplus_{\substack{I \in \II(n) }} \bigoplus_{\substack{J \in \II(n)}} W^+_{IJt} \big/ W^-_{IJt}
        \cong
        \bigoplus_{\substack{I \in \II(n) }}\bigoplus_{\substack{J \in \II(n)}} \subU_{IJt}(f).
    \]
    Now, 
    using Corollary~\ref{cor:decomposition_Z}, we obtain that if either $I \not\in\II(V)$ or $J \not\in\II(U)$, then $\subU_{IJt}(f)=0$. 
    Consequently, the first direct sum of the previous expression can be simplified to
    \[
        \bigoplus_{\substack{I \in \II(V) }}
        \bigoplus_{\substack{J \in \II(U) }}
        \subU_{IJt}(f)
        \cong 
        f V_t.
    \]
    Notice that the structure maps of $\subU_{IJt}(f)$ and $f V_t$ are both induced by the structure maps of $U$,
    and they commute with the isomorphisms given by each $t$.
    Then,
    \[
        \bigoplus_{\substack{I \in \II(V) }}
        \bigoplus_{\substack{J \in \II(U) }}
        \subU_{IJ}(f)
        \cong 
        f V.
    \]
\end{proof}

\section{The chain contraction of an edge collapse}
\label{sec:collapses}

Let $K$ be a flag complex and consider a simplicial map $p\colon L\rightarrow K$.
Let $K'$ be the complex obtained after an edge collapse of $K$, and $L'$ after one of $L$.
Then, in general, one does not need to have a simplicial map $L'\rightarrow K'$.
In particular, it is not straightforward to deduce properties about $p_*\colon \PH(L)\rightarrow \PH(K)$ from the collapsed complexes $L'$ and $K'$.
To bypass this problem, we consider the following diagram
\[
	 \begin{tikzcd}[column sep=large]
	 	C(L) \arrow[r,"p"]\arrow[d, "\pro{L}", shift left=.5ex] & C(K) \arrow[d, "\pro{K}", shift left=.5ex]\\
	 	C(L') \arrow[u, "\inc{L}", shift left=.5ex] \arrow[r, dotted 
        ] & C(K').\arrow[u, "\inc{K}", shift left=.5ex]
	 \end{tikzcd}  
 \]
 where the vertical arrows are part of chain contractions $(\pro{K}, \inc{K}, \har{K})$ and  $(\pro{L}, \inc{L}, \har{L})$ induced by performing edge collapses on both $K$ and $L$.
In particular,
$p_* \cong \left( \pro{K} \circ p \circ \inc{L} \right)_*$,
and so we can focus on understanding the chain map $\pro{K} \circ p \circ \inc{L}\colon C(L')\rightarrow C(K')$.
The upshot of this viewpoint is that one can perform independent edge collapses on the complexes $L$ and $K$ and still be able to study $p_*$ by understanding an equivalent chain map $C(L') \rightarrow C(K')$. Since, in general, the sizes of the collapsed complexes are much smaller, this presents a substantial computational boost (cf Example~\ref{ex:computation}).
\smallskip

In practice, by `studying' the persistence map $p_*$ we mean that we want to compute its associated matrix with respect to a pair of fixed persistent homology bases. 
To compute this, one varies over $d\geq 0$, and, for each representative of an element from basis for $\Ho_d(K')$, evaluate $\pro{L}_d \circ p_d \circ \inc{K}_d$ and express the image in terms of a basis of $\Ho_d(L')$.
Hence, we conclude that it is worthwhile to understand explicitly the maps composing a chain contraction $(\pro{K}, \inc{K}, \har{K})$ that results from an edge collapse operation.
Throughout the rest of this section, we focus on this question.
\smallskip

Let $e = [ u, v ]$ be an edge dominated by $v'$.
We define a degree-$1$ chain map $\har{K} \colon C(K) \to C(K) $ by setting, for a given $\sigma \in K_d$,
\[
    \har{K} (\sigma) = 
	\begin{cases}
        -[v', u,v, w_2, \ldots, w_d] & \text{ if }  
        \sigma = [u,v,w_2, \ldots, w_d]
        \text{ and } v' \notin \sigma, 
        \\
		0 & \text{ if } \{u,v\}\not\subset \sigma \mbox{ or } v' \in \sigma\ ;
	\end{cases}
\]
the definition of the chain map follows by linear extension.
In particular, note that $\har{K}_0 = 0$.
Next, we define the projection $\pro{K} \colon C(K) \to C(K')$ as follows: for  any simplex $\sigma \in K_d$, we set
	\[
        \pro{K}(\sigma) = 
		\begin{cases}
            0 & \text{if } \sigma = [ v', u, v, \ldots] \\
		    [ v', v, w_2, \ldots]  - [v', u, w_2, \ldots] & \text{if } v' \notin \sigma \; \mbox{ and }  \sigma = [u,v,w_2, \cdots, w_d]  
            \\ 
            \sigma &  
            \text{if } \{u,v\}\not\subset \sigma\ ,
		\end{cases}
	\]
    and the definition of $\pi^K$ follows by linear extension.
	In particular, notice that $\pro{K}_0 = \Id$.
    Under this definition, we claim that the triplet $(\pro{K}, \inc{K}, \har{K})$ forms a chain contraction.
\begin{lemma}\label{lem:explicit-projection}
	Let $K$ be a simplicial complex and $K'$ the resulting complex after collapsing the edge $[u,v]$.
    Then, $(\pro{K}, \inc{K}, \har{K})$ as defined above is a chain contraction.
\end{lemma}
\begin{proof}
To start, it is straightforward to check that $\pro{K}\, \har{K} = \har{K} \, \inc{K} = 0$ and $\pro{K} \, \inc{K} = \Id$.
Next, we  proceed to show that 
\begin{equation}\label{eq:is-chain-contraction}
i^K\pi^K (\sigma)=(\Id + \partial \phi^K + \phi^K\partial)(\sigma)
\end{equation}
for all simplices $\sigma \in K_d$ and all $d\geq 0$. Once we have shown that Equation~(\ref{eq:is-chain-contraction}) holds, then the result follows by linearity. 
To proceed, we check~(\ref{eq:is-chain-contraction}) in three possible cases:
\begin{enumerate}[label=\textbf{(\roman*)}]
    \item $\sigma = [v',u,v,w_3,\ldots, w_d]$. 
    In this case, by definition, $\pi^K(\sigma) = 0$ which implies that the left hand side from~(\ref{eq:is-chain-contraction}) is zero. 
    On the other hand, since $\phi^K(\sigma)=0$, we have
    \begin{align*}
		& \left( \Id + \partial \har{K} + \har{K} \partial \right) ([v', u,v,w_3, \ldots, w_{d}]) = \\
		& = \, [v', u,v,w_3, \ldots, w_{d}] + \har{K} \partial ([v', u,v,w_3, \ldots, w_{d}]).
	\end{align*}
	Now, note that
	\begin{align*}
		\partial&([v', u,v,w_3, \ldots, w_{d}]) =
		[u,v,w_3, \ldots, w_d ] - [v',v,w_3, \ldots, w_d ] \\
		&+ [v',u,w_3, \ldots, w_d ] + \sum_{i=3,\ldots,d} (-1)^i [v',u,v,w_3, \ldots, \widehat{w}_i, \ldots, w_d].
	\end{align*}
	So $\har{K} \, \partial([v', u,v,w_3, \ldots, w_{d}]) = - [v',u,v,w_3, \ldots, w_d ]$ and the right hand side from~(\ref{eq:is-chain-contraction}) also vanishes. 
    
    \item $\sigma = [u,v,w_2,\ldots, w_d]$ and $v' \not\in \sigma$.
    We observe that 
    \begin{align*}
		\partial & \har{K} ([u,v,w_2, \ldots, w_{d}]) = \,
		-\partial \left( [v', u,v,w_2, \ldots, w_{d}] \right) = \\ 
		= & \, -[u,v,w_2, \ldots, w_{d}]
		+ [v',v,w_2, \ldots, w_{d}] - [v',u,w_2, \ldots, w_{d}] \\
		& + \sum_{i=2,\ldots,d} (-1)^{i} [v',u,v, w_2, \ldots, \hat{w}_i, \ldots, w_{d}].
	\end{align*}
	and,
	\begin{align*}
		\har{K} & \partial  ([u,v,w_2, \ldots, w_{d}])  = \\
		&\, \har{K} \bigg( [v,w_2, \ldots, w_{d}] - [u,w_2, \ldots, w_{d}] 
        \\ & + \sum_{i=2,\ldots,d} (-1)^i [u,v,w_2, \ldots, \hat{w}_i, \ldots, w_{d}] \bigg) = \\
		& \, \sum_{i=2,\ldots,d} (-1)^{i+1} [v', u, v,w_2, \ldots, \hat{w}_i, \ldots, w_{d}].
	\end{align*}
	Hence, 
	\begin{align*}
		&\left( \Id  + \partial \har{K} + \har{K} \partial \right) ([u,v,w_2, \ldots, w_{d}]) = \\
 		& = \, [v',v,w_2, \ldots, w_{d}] - [v',u,w_2, \ldots, w_{d}] 
        = i^K\pro{K}([u,v,w_2, \ldots, w_{d}])\,,
	\end{align*}
    and Equation~(\ref{eq:is-chain-contraction}) also holds in this case.

    \item $\{v,u\} \not\subset \sigma$. 
    In this case, it follows that $\partial\phi^K(\sigma)= \phi^K\partial(\sigma)=0$.
    On the other hand, as $\pi^K(\sigma)=\sigma$, Equation~(\ref{eq:is-chain-contraction}) holds.  \qedhere
\end{enumerate}
\end{proof}

Notice that a sequence of chain contractions arising from collapsing edges can be concatenated to produce a valid chain contraction (see composition of chain contractions in Section~\ref{sec:chain_contraction}).
In particular, when $d=1$, $\pro{}_1$ is especially easy to evaluate: for each chain, we change the term $[u,v]$ by $[v',v] - [v',u]$.

\begin{example}
    \begin{figure}[htbp]
     \centering
     \begin{subfigure}[b]{0.22\textwidth}
         \centering
         \caption*{\large $L$}
         \includegraphics[width=\textwidth]{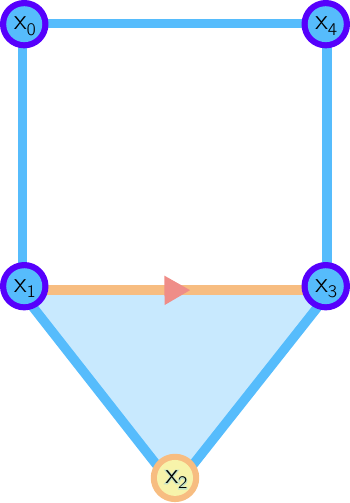}
     \end{subfigure}
     \hfill
     \begin{subfigure}[b]{0.22\textwidth}
         \centering
         \caption*{\large $L'$}
         \includegraphics[width=\textwidth]{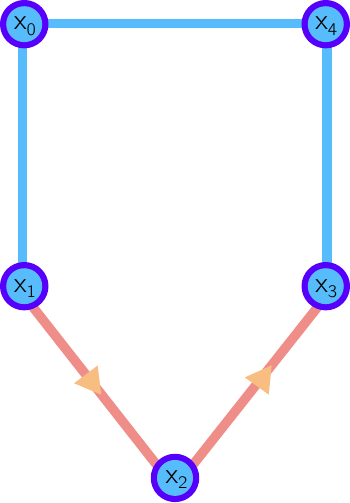}
     \end{subfigure}
     \qquad
     \begin{subfigure}[b]{0.22\textwidth}
         \centering
         \caption*{\large $K$}
         \includegraphics[width=\textwidth]{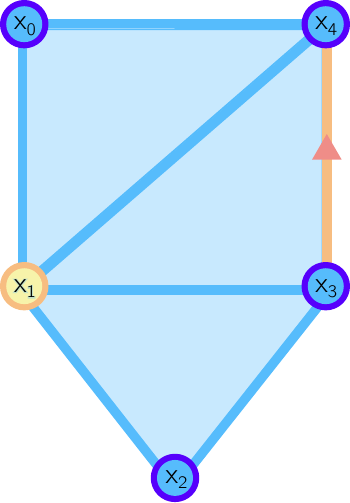}
     \end{subfigure}
     \hfill
     \begin{subfigure}[b]{0.22\textwidth}
         \centering
         \caption*{\large $K'$}
         \includegraphics[width=\textwidth]{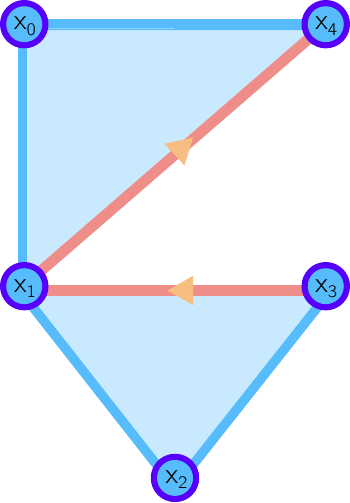}
     \end{subfigure}
        \caption{Two examples of flag complexes with a dominated edge and the corresponding complex after the collapse. The dominated edges and their projection have been highlighted.}
        \label{fig:flag_complexes}
\end{figure}

Consider the simplicial complexes from Figure~\ref{fig:flag_complexes},
and let $z$ be the $1$-chain $[ x_0, x_1] + [ x_1, x_3] +  [ x_3, x_4] - [ x_0, x_4]$ living in $C_1(K)$ and $C_1(L)$.
The edge collapse of $K$ projects $z$ to
\begin{align*}
    \pro{K}_1(z) &= [ x_0, x_1] + [ x_1, x_3] +  [ x_1, x_4] - [ x_1, x_3] - [ x_0, x_4] \\
    &= [ x_0, x_1]  +  [ x_1, x_4]  - [ x_0, x_4]
\end{align*}
in $C_1(K')$, and the one of $L$ to
\begin{align*}
    \pro{L}_1(z) &= [ x_0, x_1] + [ x_2, x_3] - [ x_2, x_1] +  [ x_3, x_4] - [ x_0, x_4] \\
    &= [ x_0, x_1] + [ x_2, x_3] + [ x_1, x_2] +  [ x_3, x_4] - [ x_0, x_4]
\end{align*} 
in $C_1(L')$.
\end{example}

\subsection{Calculating the persistence map of a flag complex}\label{sec:calc-persistence-map}
Given a pair of flag complexes, $L \subset K$, we obtain directly a persistence map between their persistent homology modules:
\[
    \begin{tikzcd}[column sep=0.75em]
        K_1 \arrow[r, phantom, "\subseteq" sloped] & K_2 \arrow[r, phantom, "\subseteq" sloped] & \ldots \arrow[r, phantom, "\subseteq" sloped] & K_n \\
        L_1 \arrow[r, phantom, "\subseteq" sloped]\arrow[u, phantom, "\subseteq" sloped] & L_2 \arrow[r, phantom, "\subseteq" sloped]\arrow[u, phantom, "\subseteq" sloped] & \ldots \arrow[r, phantom, "\subseteq" sloped] & L_n \arrow[u, phantom, "\subseteq" sloped]
    \end{tikzcd}
    \; \Rightarrow \;
    f \coloneqq
    \begin{tikzcd}[column sep=1em]
        \Ho(K_1) \arrow[r] & \Ho(K_2) \arrow[r] & \ldots \arrow[r] & \Ho(K_n) \\
        \Ho(L_1) \arrow[r]\arrow[u] & \Ho(L_2) \arrow[r]\arrow[u] & \ldots \arrow[r] & \Ho(L_n) \arrow[u]
    \end{tikzcd}
\]
On the other hand, using the edge collapse operations mentioned in Subsection~\ref{sec:flag-complexes}, we get new filtrations $K', L'$ where the addition of simplices has been delayed as much as possible (shifting) or even removed by an edge collapse in $K,L$ (trimming).
Hence, we get the same persistence module generated by smaller complexes:
\[
    \begin{tikzcd}[column sep=1em]
        K_1 \arrow[r, phantom, "\subseteq" sloped] & K_2 \arrow[r, phantom, "\subseteq" sloped] & \ldots \arrow[r, phantom, "\subseteq" sloped] & K_n \\
        K'_1 \arrow[r, phantom, "\subseteq" sloped]\arrow[u, phantom, "\subseteq" sloped] & K'_2 \arrow[r, phantom, "\subseteq" sloped]\arrow[u, phantom, "\subseteq" sloped] & \ldots \arrow[r, phantom, "\subseteq" sloped] & K'_n \arrow[u, phantom, "\subseteq" sloped]
    \end{tikzcd}
    \quad \Rightarrow \quad
    \begin{tikzcd}[column sep=1em]
        \Ho(K_1) \arrow[r]\arrow[d, phantom, "\cong" sloped] & \Ho(K_2) \arrow[r]\arrow[d, phantom, "\cong" sloped] & \ldots \arrow[r] & \Ho(K_n)\arrow[d, phantom, "\cong" sloped] \\
        \Ho(K'_1) \arrow[r] & \Ho(K'_2) \arrow[r]& \ldots \arrow[r] & \Ho(K'_n)
    \end{tikzcd}
\]
In particular, every $K'_i$ is obtained from $K_i$ after a sequence of edge collapses.
Denote by $(\pro{K_i}, \inc{K_i}, \har{K_i})$
the chain contraction obtained from these sequences of edge collapses for each $i$, and $p^i$ to the inclusion of $L_i$ into $K_i$.
Then,
\[
    \begin{tikzcd}
	   U \\
	   V \arrow[u, "f"]
    \end{tikzcd} \cong
    \begin{tikzcd}[column sep=1em]
        \Ho_d(K'_1) \arrow[r] & \Ho_d(K'_2) \arrow[r] & \ldots \arrow[r] & \Ho_d(K'_n) \\
        \Ho_d(L'_1) \arrow[r]\arrow[u, "(\pro{L_1}_d \, p^1_d \, \inc{K_1}_d)_*", swap] & \Ho_d(L'_2) \arrow[r]\arrow[u, "(\pro{L_2}_d \, p^2_d \, \inc{K_2}_d)_*", swap] & \ldots \arrow[r] & \Ho_d(L'_n) \arrow[u, "(\pro{L_n}_d \, p^n_d \, \inc{K_n}_d)_*", swap]
    \end{tikzcd}
\]
To obtain the matrix representation, $F$, directly from this diagram, we just have to get bases $\cA = [ \alpha_\imath]$ and $\cB = [ \beta_\jmath]$ for $V$ and $U$, and obtain the indices as in Definition~\ref{def:assoc-matrix}:
\[
    f\big(\alpha_{\jmath a_\jmath}^1\big) = \big(
    \pro{L_{a_\jmath}} \, p^{a_\jmath} \, \inc{K_{a_\jmath}}
    \big)_*\big(\alpha_{ja_\jmath}^1\big) = \sum_{\beta_{\imath a_\jmath} \in \cB_{a_\jmath}} F(i,j) \beta_{\imath a_\jmath}^1.
\]

\begin{example}
    Consider the flag complexes appearing in Figure~\ref{fig:flag_complexes}.
    Note that $L$ is contained in $K$ via an inclusion $p$.
    However, $L'$ is no longer contained in $K'$.
    Consider now the $1$-chain $z' = [ x_0, x_1] + [ x_1, x_2] + [x_2, x_3] + [ x_3, x_4] - [ x_0, x_4]$ in $C_1(L')$. 
    What is the representative of $z'$ in $C_1(K')$?
    It is given by
    \begin{align*}
        \pro{K}_1 p \, \inc{L}_1 (z')
        =&
        [ x_0, x_1] + [ x_1, x_2] + [ x_2, x_3] -  [ x_1, x_3] + [ x_1, x_4]  - [x_0, x_4].   
    \end{align*}
\end{example}

\subsection{Cycle collapse algorithm}

In this subsection, we present a procedure that collapses representative cycles following shifting or trimming operations; c.f. Subsection~\ref{sec:flag-complexes}.
In particular, this procedure is meant as a subroutine of the ``core flag filtration backward algorithm'' which is Algorithm~1 from~\cite{swap}.
Basically, we consider a list of representative cycles and update them each time an edge is either shifted or trimmed.
This is summarized in Algorithm~\ref{alg:PH-collapse} which is based on the definition of $\pro{K}$; being the projection of a chain contraction by Lemma~\ref{lem:explicit-projection}. 
\smallskip

As before, let $L$ and $K$ be two flag complexes together with a simplicial map $p\colon L\rightarrow K$ and let $L'$ be the complex obtained after an edge collapse of $L$. 
The method described here (summarized in Algorithm~\ref{alg:PH-collapse}) is applied immediately prior to collapsing the edges of $K$.
In particular, suppose that persistent homology $\PH(L')$ has been computed, leading to barcodes and representative cycles. 
We denote by $\repsL$ the representative cycles obtained by means of this procedure. 
That is, $\repsL$ is a list where each element $\gamma \in \repsL$ is a chain in $C_d(L)$ representing some persistent homology class.
We also denote by $t(\gamma)$ the birth value of $\gamma$ in $L'$, 
i.e. $t(\gamma)$ is the smallest filtration value $t$ such that $\gamma \in L'_t$.
\smallskip

Next, since $p\colon L\rightarrow K$ is a simplicial map, we compute $\imReps$; which is the set $p(\repsL)=\{p(\gamma) \mid \gamma \in \repsL\}$. 
Then, as the complex $K$ is collapsed into $K'$ via Algorithm~1 from~\cite{swap}, we proceed to collapse the cycles from $\imReps$ simultaneously. 
More concretely, each time an edge $e=\{u,v\}$ is trimmed or shifted due to a dominating vertex $v'$, we collapse the representatives from $\imReps$ via  $\pro{K}$. 
This procedure is summarized in Algorithm~\ref{alg:PH-collapse} which is explained next.
\smallskip

Let $\prevtime$ be the filtration value of $e$ before being shifted/trimmed to a filtration value $\newtime$, possibly $\infty$.
Then, Algorithm~\ref{alg:PH-collapse} proceeds as follows: for each $\gamma \in \imReps$, if $t(\gamma) \in [\prevtime, \newtime)$ then $\gamma$ is collapsed to $\pro{K}(\gamma)$, which is computed using linearity of $\pi^K$, as done in lines \ref{if:value-gamma}-\ref{get:newgamma}. 
Otherwise, if $t(\gamma) \notin [\prevtime, \newtime)$, we leave $\gamma$ unchanged. 
The resulting new representatives are returned in a list named $\repsCollapsed$.
Finally, we explain the notation used in lines~\ref{for:sigma-gamma}-\ref{gets:collapse}: using the standard basis, we write $\gamma = \sum_{\tau \in K} \gamma_\tau \tau$ for coefficeints $\gamma_\tau \in k$ for all $\tau \in K$. Then, we write $\sigma \in \gamma$ if and only if $\gamma_\sigma\neq 0$.

\begin{algorithm}[hbt]
  \caption{Cycle collapse}\label{alg:PH-collapse}
  \begin{algorithmic}[1]
    \STATE \textbf{Input:  \imReps, $\{u,v\}$, $v'$, \prevtime, \newtime}
    \STATE \textbf{Output:} Collapsed cycle representatives $\repsCollapsed$.
    \STATE {$\repsCollapsed \gets \{\}$}
	\FOR {$\gamma \in \imReps$} 
            \IF {\label{if:value-gamma} $t(\gamma) \in [\prevtime, \newtime)$}
                \STATE $\gammanew \gets 0$
                \FOR {$\sigma \in \gamma$ \label{for:sigma-gamma}} 
                    \STATE \label{gets:collapse} $\gammanew \gets \gammanew + \gamma_\sigma \pi^K(\sigma)$
                \ENDFOR
                \STATE { \label{get:newgamma}\repsCollapsed.add($\gammanew$)}
            \ELSE
            \STATE {\repsCollapsed.add($\gamma$)}
            \ENDIF
	\ENDFOR
    \STATE \textbf{Returns:} $\repsCollapsed$
  \end{algorithmic}
\end{algorithm}

\begin{remark}
    The edge collapse procedure described in \cite{swap} considers a sequence of edges which respects filtration values. In this setup, a parallelization procedure is proposed: first, the sequence of edges is divided into two halves, second, a simultaneous collapse on both halves (using two processors) is performed, and, third, a final collapse that merges both results is executed. 
    This procedure can be applied to our proposed method for collapsing representatives, as far as the last edge from the first sequence and the first edge from the second do not share the same filtration value.
    Indeed, the conditional from line~\ref{if:value-gamma} in Algorithm~\ref{alg:PH-collapse} ensures that each representative cycle $\gamma\in\imReps$ will have filtration value $t(\gamma)$ in only one of the two divisions; hence, it is collapsed by a single processor at a time.
\end{remark}

\begin{remark}
    Following the discussion on parallelisation, notice that the for loop over all $\imReps$ in line~3 from Algorithm~\ref{alg:PH-collapse} can be run in parallel.
\end{remark}

\section{Implementation and experiments}
\label{sec:experiments}

In this Section, we provide 
details about the implementation of the computation of block functions, as well as the computation of the associated matrix of a persistence map via edge collapses. 
Also, we present some experiments illustrating induced matchings.

\subsection{Overall procedure for flag complex filtrations}

Here we briefly describe the procedure employed to compute the induced matching in the case of two filtered flag complexes, $L$ and $K$, together with a simplicial map $p\colon L\rightarrow K$.
It turns out that this procedure is complex and requires of several subroutines, which we briefly describe next:
\begin{itemize}
    \item $\collapseEdges(L_1)$: given a one skeleton filtration $L_1$, it returns its core $L_1'$ using the ``core backward filtration algorithm'', see Algorithm~1 from \cite{swap}.
    \item $\flagComplex(L_1')$: given a one skeleton filtration $L_1'$, it returns the corresponding flag complex filtration $L'$. 
    \item $\computePersistence(L')$: Computes (and returns) the barcodes $\barcodeL$ and representatives $\repsL$ for the persistent homology module $\PH(L')$. 
    As a byproduct, since we are using PHAT, we also obtain the reduced differential matrix, which we denote by $R_L$.
    \item $p(\repsL)$: image under $p\colon L\rightarrow K$ for all elements from $\repsL$.
    \item $\collapseEdges(K_1, \imReps)$: given a one skeleton filtration $K_1$ and a list of cycles $\imReps$; it returns the core $K'_1$ of $K_1$ according to Algorithm~1 from \cite{swap}. 
    As it performs edge collapses, Algorithm~2 from this text is used to update the representatives from $\imReps$ accordingly. The resulting representatives are returned as a list called $\repsCollapsed$.

    \item $\getPmMatrix(R_K, \repsK, \repsCollapsed)$: returns a matrix $F$ associated to the persistence map $\PH(L)\rightarrow \PH(K)$, c.f. Definition~\ref{def:assoc-matrix}. 
    In a nutshell, we use PHAT to express $\repsCollapsed$ as a combination of $\repsK$ and $R_K$. This allows to obtain the coefficients of the homology class of each element from $\repsCollapsed$ with respect to the computed homology classes from $\PH(K)$.

    \item $\columnReduction(F)$: reduces the matrix $F$ using Algorithm~\ref{alg:column_reduction}. 

    \item $\getMatching(R, \barcodeL, \barcodeK)$: uses the reduced matrix to obtain the matching $\matching{}$ by using eq.~(\ref{eq:definition}).
    In particular, $\barcodeL$ and $\barcodeK$ are associated to the columns and rows respectively.
\end{itemize}
These subroutines are combined together into Algorithm~\ref{alg:overall}.

\begin{algorithm}[hbt]
  \caption{Induced matching}\label{alg:overall}
  \begin{algorithmic}[1]
    \STATE \textbf{Input:}  $L_1,K_1$ two one skeleton filtrations and $p\colon L_1\rightarrow K_1$.
    \STATE \textbf{Output:} barcodes and matching of  $p\colon \PH(L)\rightarrow \PH(K)$.
    \STATE {$L_1' \gets \collapseEdges(L_1)$}
    \STATE {$L' \gets \flagComplex(L_1')$}
    \STATE {$\barcodeL, \repsL, R_L\gets \computePersistence(L')$}
    \STATE {$\imReps \gets p(\repsL)$}
    \STATE {$K_1', \repsCollapsed \gets \collapseEdges(K, \repsL)$}
    \STATE {$K' \gets \flagComplex(K_1')$}
    \STATE {$\barcodeK, \repsK, R_K \gets \computePersistence(K')$}
    \STATE {$F \gets \getPmMatrix(R_K,\repsK, \repsCollapsed)$}
    \STATE {$R \gets \columnReduction(F)$}
    \STATE {$\matching{} \gets \getMatching(R, \barcodeL, \barcodeK)$}
    \STATE {\textbf{Returns:} $\barcodeL, \barcodeK, \matching{}$}
  \end{algorithmic}
\end{algorithm}

\subsection{Implementation details for Vietoris-Rips filtrations}

We have implemented an adaptation of the edge collapse module from Gudhi, that we have called ``permovec''\footnote{https://bitbucket.org/atorras1618/permovec}, standing for \textbf{per}sistence \textbf{mo}rphisms \textbf{v}ia \textbf{e}dge \textbf{c}ollapses.  
The procedure is briefly described next. 
First, one computes the Vietoris-Rips complex and its corresponding edge collapse by using the respective modules from \cite{gudhi}. 
Next, given the one skeleton of the collapsed complex, we expand its simplices up to dimension two (this is done by using the simplex tree module from Gudhi). 
Then, this object is used to compute its boundary matrix, which is stored as an object from the library  PHAT: \cite{phat}. 
This allows the recovery of representative cycles for the persistent homology classes in the domain. Such cycles are then collapsed by using the transformations described in Section~\ref{sec:collapses}, and, in particular, Algorithm~\ref{alg:PH-collapse}. 
\smallskip

Following this step, the matrix associated with a persistence map is recovered by performing column reductions on PHAT boundary matrix objects.
For practical reasons, the current implementation is restricted to dimensions $0$ and $1$ with field coefficients in $\mathbb{Z}_2$ (see Remark~\ref{rem:phat}).
Future updates will extend the code to higher dimensions and arbitrary field coefficients, aligning with the full generality of Algorithm~\ref{alg:PH-collapse}.

\smallskip

\begin{remark}\label{rem:phat}
We have slightly modified the PHAT library so that it allows to track the performed column additions when reducing, and, in particular, the resulting representative cycles\footnote{https://bitbucket.org/atorras1618/phat/src/master}.
Our choice of using the PHAT library forced us to implement the procedure for $\mathbb{Z}_2$ coefficients, since the PHAT library only supports this coefficient field. 
\end{remark}
\smallskip

The computation of the induced matching is then implemented via Python code that is easy to use, called ``IBloFunMatch''\footnote{https://github.com/Cimagroup/IBloFunMatch}, standing for induced block function and matchings. 
In this repository the user can find Python code for easily using permovec and computing also the induced matching described in Section~\ref{sec:induced}. In particular, all experiments from this article can be found within a folder named ``Notebooks-JoSC''.
An aspect that is considered in this implementation is that, as explained in~\cite{edge_collapse}, a single iteration over all edges might not exhaust all possible collapses.
This is why the user can specify the number of iterations over all edges to perform when collapsing the edges.

\subsection{Experiments}

This section begins with an example of induced matchings. 
We then show the computational advantages of employing edge collapses compared to standard approaches.

\begin{example}
    \begin{figure}
    \centering
    \includegraphics[width=0.5\linewidth]{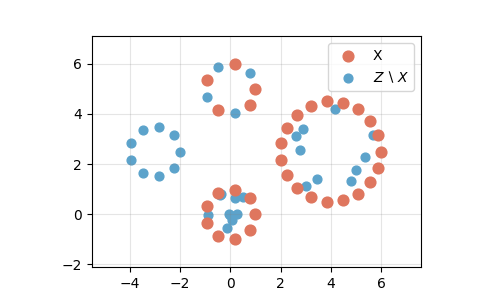}
    \caption{Depiction of a pair of point samples $X \subset Z$ where the points from $X$ are plotted in red while points from $Z \setminus X$ are plotted in blue.}
    \label{fig:points-circles}
    \end{figure}
    Consider the pair of point samples $X\subset Z$ depicted in Figure~\ref{fig:points-circles}. 
    As described in Section~\ref{sec:collapses}, one starts by performing edge collapses on $L=\VR(X)$ and then computing the persistent homology of $L'$. In this case, there are three cycles in dimension $1$ whose representatives are depicted on the left of Figure~\ref{fig:matching-circles-reps}. 
    Next, one proceeds to collapse the simplicial complex $K=\VR(Z)$ by performing edge collapses and, in doing so, also collapsing the representative cycles from $L'$; as described in Section~\ref{sec:calc-persistence-map}. 
    These representatives are then projected into $K'$ via composing the respective maps (one for each edge collapse) $p_*$ described in Section~\ref{sec:collapses}.
    This allows to obtain the matrix $F$ associated to $\PH_1(L')\rightarrow \PH_1(K')$ depicted on the left of Figure~\ref{fig:matching-circles}.
    In this case, the matching can be read directly from $F$ since pivots do not repeat and the resulting induced matching $\matching{f}$ is depicted on the right of Figure~\ref{fig:matching-circles}.
    In addition, we also depict the matched representative cycles on the right of Figure~\ref{fig:matching-circles-reps}.
    \exend
    \begin{figure}
    \centering
    \includegraphics[width=0.8\linewidth]{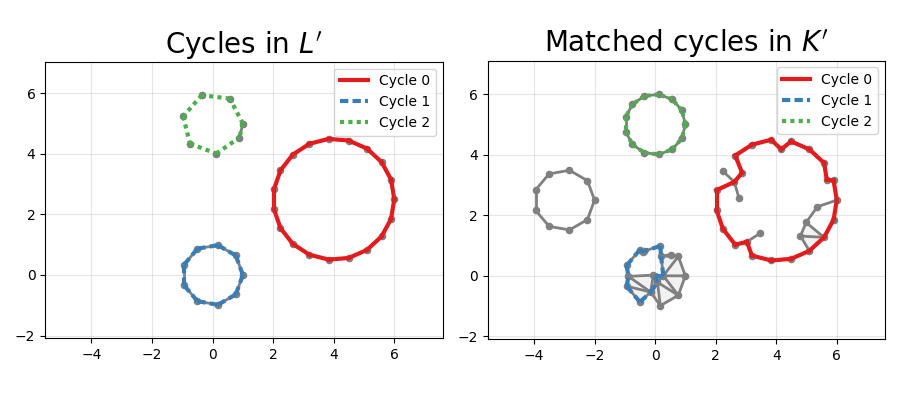}
    \caption{Depiction of three cycles in $L'$ and in $K'$ corresponding to the matching $\matching{f}$.
    The flag complexes $L'$ and $K'$ are plotted in gray up to a fixed filtration value.}
    \label{fig:matching-circles-reps}
    \end{figure}
    \begin{figure}
    \centering
    \begin{subfigure}{0.15\textwidth}
    {\scriptsize
        \begin{tikzpicture} 
        \node at (0,0) {$
            F=\left[\begin{array}{ccc}
                0 & 0 & 0 \\
                0 & 0 & 0 \\
                0 & \textbf{1} & 0 \\
                0 & 0 & 0 \\
                0 & 0 & 0 \\
                0 & \textbf{1} & 0 \\
                0 & 0 & 0 \\
                0 & 0 & \textbf{1} \\
                \textbf{1} & 0 & 0 \\
                0 & 0 & 0 
            \end{array}
            \right]
            $};
        \end{tikzpicture}
    }
    \end{subfigure}
    \hfill
    \begin{subfigure}{0.75\textwidth}
        \includegraphics[width=0.95\linewidth]{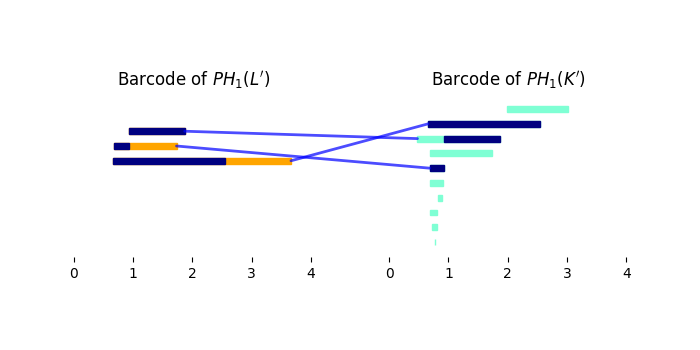}
    \end{subfigure}
    \caption{Matrix associated to $\PH_1(L')\rightarrow \PH_1(K')$ (left) and corresponding induced matching (right). 
    The order of columns (from left-to-right) and rows (from top-to-bottom) in $F$  corresponds to the order of plotted intervals on the right (from bottom-to-top).
    In addition, we mark (in dark blue) the intersection between matched intervals.}
    \label{fig:matching-circles}
    \end{figure}
\end{example}

\begin{example}\label{ex:computation}
    \begin{figure}
        \centering
        \includegraphics[width=0.3\linewidth]{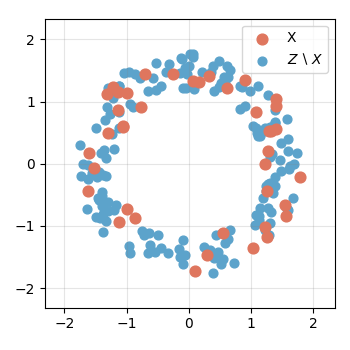}
        \caption{Pair $X \subset Z$ of sampled points around a circle.}
        \label{fig:example-computation-points}
    \end{figure}
    We consider a point sample $Z$ of $300$ points distributed around a circle and a subset $X\subset Z$ as depicted in Figure~\ref{fig:example-computation-points}.
    In this example, we study the advantage of using edge collapses by comparing computations performed on a standard desktop computer
    \footnote{We performed the computations on a Host OS Microsoft Windows 11 Home, with 32 GB RAM and a 12th Gen Intel(R) Core(TM) i7-12700H (2.30 GHz) processor. 
    The computations were run on WSL2. }. 
    For example, in the top of Figure~\ref{fig:cycles-example-computation} we observe the representative cycles for $L'$, their image in $K'$ (after projecting using the edge collapses) and their matched representatives. 
    Executing the edge collapses, calculating the resulting persistence map matrix and the induced matching took around $0.63$ seconds. 
    On the bottom of Figure~\ref{fig:cycles-example-computation} we depict the representative cycles for homology classes for $L$, their embedding into $K$ and their matching representatives in $K$, without using edge collapses. 
    Calculating the persistence map  without the edge collapses and the induced partial matching took about 29 seconds. 
    This represents a significant speedup, confirming the efficiency gains of the edge collapse approach. 
    Also, edge collapses lead to a reduced memory usage; e.g. for a 180 point sample the peak memory used is reduced from 770 Mb down to 3 Mb with a single edge collapse iteration. For measuring memory we used the ``massif'' tool from Valgrind~\citep{Valgrind2006}.
    All these performance measures are summarized in Figure~\ref{fig:computation-comparison}. 
    \exend
    \begin{figure}
        \begin{subfigure}{\textwidth}
        \centering
            \includegraphics[width=0.9\linewidth]{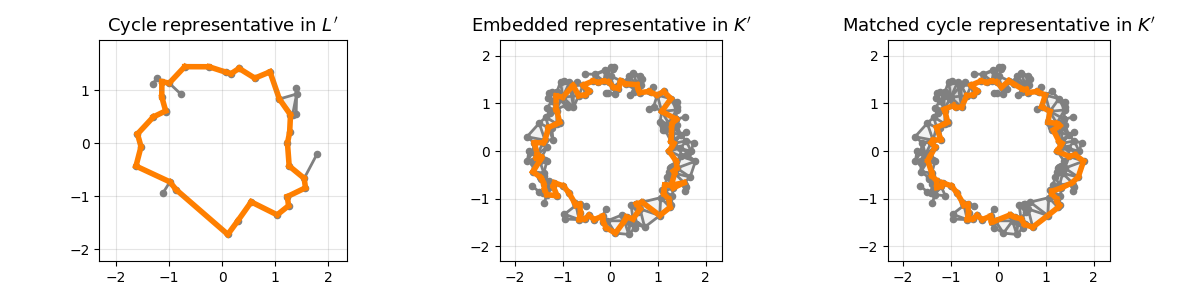}
        \end{subfigure}
        \begin{subfigure}{\textwidth}
        \centering
            \includegraphics[width=0.9\linewidth]{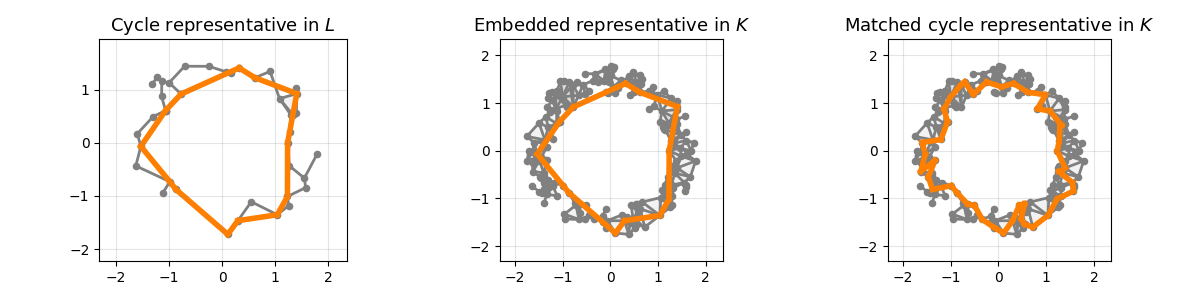}
        \end{subfigure}
        \caption{
        Depiction of representative cycles for the same homology class in $L'$ and $L$ (left), its projection or embedding into $K'$ and $K$ (middle) and the cycle it is matched to (right).
        In particular, the top row illustrates the matching of representative cycles using edge collapses, while, in the bottom row, collapses were not used.}
        \label{fig:cycles-example-computation}
    \end{figure}
    \begin{figure}
        \centering
        \includegraphics[width=0.49\linewidth]{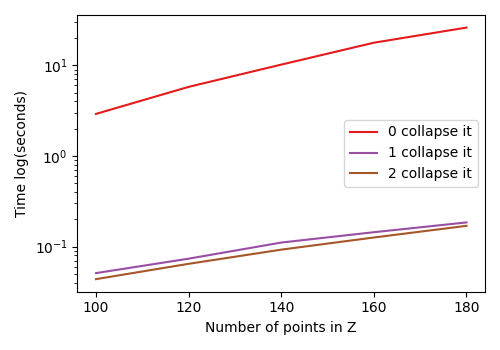}
        \includegraphics[width=0.49\linewidth]{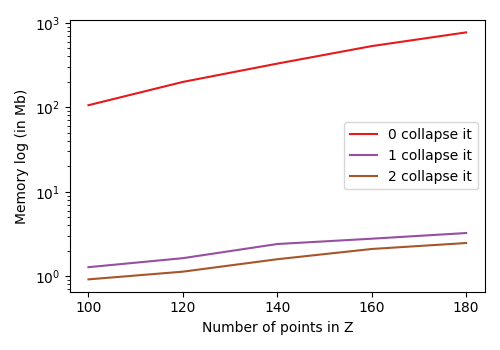}
        \caption{Runtimes (left) and peak memory usage (right) for computing the associated matrix of a persistence map and the induced matching. In logarithmic scale, a significant improvement can be noticed when a single edge collapse iteration is performed (over not collapsing edges). An additional edge iteration seems to also improve performance slightly. }
        \label{fig:computation-comparison}
    \end{figure}
\end{example}

\section{Conclusions and Future work}
We have introduced a matrix reduction algorithm that takes a persistence map in matrix form and calculates a partial matching between the persistence diagrams.
This partial matching also provides the image module, making it a richer invariant, and can be seen as a projection of persistence maps to interval decomposable maps.
In addition, it can be calculated in cubic time and is additive with respect to direct sums.
\smallskip

For persistent homology of flag complexes, we provide a chain contraction that allows for the calculation of persistence maps directly from edge-collapsed flag complexes, considerably speeding up the calculations.
We have also added computational examples illustrating calculation of the persistence map combined with the induced partial matching. 
\smallskip

There are two appendices where we explain the differences of our method with $\chi_f$, the most common induced partial matching, and $\cM_f$, a block function that inspired our method.
In particular, we explain the difference between the methods and provide a small space where all methods coincide.

A question left to future research is whether one can obtain an approximate matching in an analogous way to the approximate barcode obtained via edge collapses in \cite{swap}. 
We also plan to investigate possible applications of $\matching{f}$ in topological data analysis and its relation with machine learning.
Initial progress in data quality \citep{quality} and agent fleets \citep{europe} provides a basis for this work.

As commented in the introduction, there are many methods in the literature with a similar aim but different in spirit.
As future work, it would be interesting to study how these invariants related to each other, a nontrivial task due to the variety of concepts used to define each of them.
\smallskip

In addition, we believe that this method can be extended to other cases, like two persistence maps with a common domain: $V \leftarrow W \rightarrow U$.
These persistence maps appear in dynamical systems \citep{self, sampled} and can be used in a broader concept than persistence maps \citep{extension}.

\section*{Acknowledgements}
This project was partially funded by the  MCIN/AEI and the NextGenerationEU/PRTR, under project TED2021-129438B-I00.
The authors thank IMUS-Maria de Maeztu grant CEX2024-001517-M - Apoyo a Unidades de Excelencia María de Maeztu for supporting this research, funded by MICIU/AEI/ 10.13039/501100011033.
The authors would also like to thank Lars M Salbu for fruitful discussions regarding the operators from Definition~\ref{def:V-tilde-hat} and their relation with the order relations introduced in Definition~\ref{def:orders}.

\bibliographystyle{elsarticle-harv} 
\bibliography{references}

\appendix
\section{Comparison between $\matching{f}$ and $\chi_f$}
\label{app:BL}

In this appendix we compare $\matching{f}$ with the most common induced partial matching in the literature, $\chi_f$.
It was first introduced to simplify the proof of the stability theorem~\citep{BauLes2015, BauLes2020}, but has also been used in TDA applications~\citep{faithful, ReBo23}.
In this section, we will see that both partial matchings have a different nature, but there is a family of persistence modules for which they coincide.
\smallskip

Before continuing, we recall how to calculate $\chi_f$, which is defined on sets instead of multisets.
The set representation of the persistence diagram $\PD(V)$ is denoted as \emph{barcode} of $V$, $\mtb(V)$.
\[
     \mtb(V)
     \coloneqq
     \{I_i \mid I \in \II(V), i \in [m_s] \}.
\]
Given $e\in [n]$, we say that $I_i \in \mtb(V)^e \subset \mtb(V)$ if and only if $e$ is the right endpoint of the interval $I_i$. 
Analogously, 
$I_i\in \mtb(V)_e \subset \mtb(V)$ if and only if $e$ is the left endpoint of $I$.
Observe that we can fix an order in each of these sets, for $\mtb(V)^e$, 
$[a,e]_i <^e [b,e]_j$ if $a < b$ or $a = b$ and $i < j$; and for $\mtb(V)_e$,
$[e,a]_i <_e [e,b]_j$ if $a > b$ or $a = b$ and $i < j$.
The definition of $\chi_f$ is based on the following result.
 
\begin{theorem}[Theorem 4.2 from~\citep{BauLes2015}]\label{th:bauer}
    If $f : V \rightarrow U$ is injective, then for each $e \in [n]$ we have that
    \[
        \#\mtb(V)^e \leq \# \mtb(U)^e
    \]
    and if $f : V \rightarrow U$ is surjective, then
    \[
        \#\mtb(V)_e \geq \# \mtb(U)_e.
    \]
\end{theorem}

Now, observe that, given $f\colon V \rightarrow U$, there exists a unique decomposition $V \overset{f'}{\twoheadrightarrow} f V \overset{f''}{\hookrightarrow} U $ where $f'$ is surjective and $f''$ is injective~\citep{BauLes2015, BauLes2020}.  
Applying Theorem~\ref{th:bauer},
we have that $\#\mtb(V)_e \geq \#\mtb(f V)_e$ for each $e \in [n]$.
Using the defined order for these sets, 
we can build a partial matching $\sigma_1: \mtb(V) \to \mtb(fV)$ that maps
the $j$-th element of $\mtb(V)_e$ to the $j$-th element of $\mtb(f V)_e$ for all $e \in [n]$.
This way,
all the bars in $\mtb(f V)$ are matched.
Analogously, we can match the $j$-th element of $\mtb(f V)^e$ to the $j$-th element of $\mtb(U)^e$,
obtaining an injection $\sigma_2:\mtb(fV)\to \mtb(U)$.
Then, a partial matching $\sigma$ between $\mtb(V)$ and $\mtb(U)$ can be built composing $\sigma_1$ and $\sigma_2$.

\begin{example}\label{ex:calc_chi}
    Going back to Example~\ref{exa:limitation}, recall
     \begin{equation*}		
    \begin{tikzcd}
	U \\
	V \arrow[u, "f"]
    \end{tikzcd}
    \cong 
    \begin{tikzcd}
        k \arrow[r, "\Id"] & k \arrow[r] & 0 \\
        0 \arrow[r]\arrow[u] & 0 \arrow[r]\arrow[u] & 0\arrow[u]
    \end{tikzcd}
    \oplus
    \begin{tikzcd}
        0 \arrow[r] & 0 \arrow[r] & 0 \\
        0 \arrow[r]\arrow[u] & k \arrow[r, "\Id"]\arrow[u] & k \arrow[u]
    \end{tikzcd}
    \oplus
    \begin{tikzcd}[/tikz/column 3/.style={column sep=-0.5em}]
        k \arrow[r, "\Id"] & k \arrow[r] & 0 \\
        0 \arrow[r]\arrow[u] & k \arrow[r]\arrow[u, "\Id"] & 0\arrow[u] & .
    \end{tikzcd}
    \end{equation*}
    Note that $\mtb(V) = \{[2,2]_1, [2,3]_1\}$,  $\mtb(U) = \{[1,2]_1, [1,2]_2\}$, and $\mtb(f V) = \{[2,2]_1\}$.
    In addition, the order in $\mtb(V)_2$ establishes that $[2,3]_1 <_2 [2,2]_1$,
    and the one in $\mtb(U)^2$ establishes that $[1,2]_1 <^2 [1,2]_2$.
    Then, the partial matching $\chi_f$ is given by
    \begin{equation*}
        [2,3]_1 \mapsto [1,2]_1 , \quad \emptyset \mapsto [1,2]_2 \quad \text{and} \quad [2,2]_1 \mapsto \emptyset. \exend
    \end{equation*}
\end{example}

As explained in Section~\ref{sec:introduction}, $\chi_f$ may produce outputs that do not align with the decomposition of the persistence map $f$.
However, we could construct an additive operator $\chi'_f$ considering the unique decomposition of $f$ into 
indecomposable modules
$f_1 \oplus \ldots \oplus f_n$, and then combining each $\chi_{f_i}$ to obtain the partial matching.
The main obstacle to this construction is that, as mentioned in the introduction, there is currently no efficient way to compute 
indecomposable modules.
Since we have provided a matrix algorithm to calculate $\matching{f}$,  it is natural to ask if $\matching{f}$ is equivalent to $\chi'_f$, as this would offer a way to compute an additive version of $\chi_f$.
As this example shows, this is not the case.

\begin{example}\label{ex:m_isnot_chi}
Consider the following 
indecomposable module:
    \[\
        \begin{tikzcd}
    	U \\
    	V \arrow[u, "f"]
        \end{tikzcd}
        \cong 
        \begin{tikzcd}
        k \arrow[r, "\sbm{0 \\ 1}"] & 
        k^2 \arrow[r, "\sbm{0 \, 1}"] & 
        k \arrow[r] & 
        0 
        \\
        0 \arrow[r]\arrow[u] & 
        k^2 \arrow[r, "\Id"]\arrow[u,  "\sbm{1 \,0 \\ 1 \, 1}"] & 
        k^2 \arrow[r, "\sbm{0 \, 1}" 
        ]\arrow[u, "\sbm{1 \, 1}"] & 
        k\arrow[u].
    \end{tikzcd}
    \] 
    We have that $\mtb(V)=\{[2,3],[2,4]\}$, $\mtb(U)=\{[2,2],[1,3]\}$ and
    $\mtb(f V) = \{ [2,2]$, $[2,3] \}$.
    Hence, the partial matching given by $\chi_f$ is
    \[
        [2,3] \mapsto [2,2] \quad \text{ and } \quad [2,4] \mapsto [1,3].
    \]
    However, if we order the generators of $V$ as $\alpha_1\sim [2,3]$ and $\alpha_2\sim [2,4]$, and the ones of $U$ as $\beta_1\sim [2,2]$ and $\beta_2 \sim [2,3]$. 
    The matrix $F$ associated with $f\colon V\rightarrow U$ on this choice of persistence bases as well as its Gaussian reduction $R$ are
    \[
    F= 
    \left( \begin{array}{cc}1 & 0 \\ 1 & 1 \end{array}\right)\,
    \longrightarrow 
    R = \left( \begin{array}{cc}1 & 1 \\ 1 & 0 \end{array}\right).
    \]
    and the partial matching given by $\matching{f}$ is
    \[
        [2,3] \mapsto [1,3] \quad \text{ and } \quad [2,4] \mapsto [2,2].
    \]
\end{example}

Hence, both partial matchings differ from each other, and $\matching{f}$ is not the additive version of $\chi_f$. 
However, there exists a family of persistence maps for which both matchings coincide.
A persistence map $f \colon V \to U$ is called polybar if all intervals from $\PD(V)$ have multiplicity one and distinct left endpoints, while all intervals from $\PD(U)$ have multiplicity one and distinct right endpoints.
In particular, the definition of polybar map guarantees that every element of $\PD(f V)$ has multiplicity one, and there is a unique element with the same left endpoint in $\PD(V)$, and a unique element with the same right endpoint in $\PD(U)$.
\begin{theorem}
    If $f$ is a polybar map, $f V$, $\matching{f}(I,J)$ and $\chi_f$ are equivalent invariants.
    In particular, the following three statements are equivalent:
    \begin{itemize}
        \item $I \cap J \in \II(f V)$,
        \item $\chi_f(I) = J$, and
        \item $\matching{f}(I,J) = 1$.
    \end{itemize}
\end{theorem}
\begin{proof}
    Assume that $I \cap J \in \II(f V)$ implies $\chi_f(I) = J$.
    Note that, by definition of $\chi_f$, there is a pair of matched bar for every element in $\mtb(f V)$.
    Since $f$ is polybar, there is no interval repeated in $\mtb(f V)$, so it is bijective to $\II(f V)$.
    Moreover, for every interval $K \in \II(f V)$ there is a unique interval $I \in \PD(V)$ with the same left endpoint, and a unique interval $J \in \PD(U)$ with the same right endpoint.
    Then, necessarily, $I \cap J = K$ and $\chi_f(I) = J$.
    In addition, by Theorem~\ref{the:image}, $\matching{f}(I,J) = 1$.
    By the same argument, if $I \cap J \notin \II(f V)$, necessarily $\chi_f(I) \neq J$ and $\matching{f}(I,J) = 0$.
\end{proof}

We finish the section noting that $\chi_f$ is additive when $f$ is polybar.
This is not a difficult fact yet we could not find a proof in the literature.
Given two partial matchings defined on barcodes as $\xi_f$, $\sigma_1: B_1 \rightarrow B_1'$ and $\sigma_2: B_2 \rightarrow B_2'$, we denote by $\sigma_1 \sqcup \sigma_2 : B_1 \sqcup B_2 \rightarrow B'_1 \sqcup B_2'$ the new partial matching given by the disjoint union of the barcodes $B_1$ and $B_2$, and their relations.
\begin{corollary}
    If $f : V \rightarrow U$ is a polybar map and $f \cong f_1 \oplus f_2$, then $\chi_{f} = \chi_{f_1} \sqcup \chi_{f_2}$.
\end{corollary}
\begin{proof}
    This is a direct consequence of Theorem~A.4, since the image module of $fV$ is $f_1 V_1 \oplus f_2 V_2$.
\end{proof}

From this section, two interesting open questions arise:
(1) Is it possible to define a new induced matching, similarly to $\matching{f}$, that is an additive version of $\chi_f$? 
Note also that $f V$ can be obtained from $\matching{f}$, and $\chi_f$ can be obtained from $f V$. 
Thus, $\chi_f$ is completely determined by $\matching{f}$. (2) Can this fact be used to derive new properties of the relation between $\chi_f$ and $\matching{f}$?

\section{Comparison between $\matching{f}$ and $\mathcal{M}_f$}
\label{app:BF}
In~\citep{InducedMatchings2022},
the authors introduced \emph{block functions}, a weaker version of \emph{partial matchings}.
They are defined as in Definition~\ref{def:partial_matching}, but dropping inequality (\ref{eq:ine2}). Intuitively, a partial matching represents a bijection defined on two subsets, while a block function represents a function defined on two subsets. 
More precisely,
the induced block function $\cM_f$ is built similarly to $\matching{f}$ but using persistence modules different from $\subU_{IJ}$ \citep{InducedMatchings2022}.
Actually, it can be seen as the predecessor of $\matching{f}$
\smallskip

Recall from Section~\ref{sec:properties} that, given a persistence basis $\cA$ of $U$, we can define the
\[
    \cA_{I}\coloneqq \big\{ 
        \alpha_i \in \cA \mid \alpha_i \sim I
    \big\}.
\]
We also recall from Section~\ref{sec:persistence-modules} that $J=[c,d] \leq [a,b]=I$ if $c \leq a \leq d \leq b$.
To calculate $\cM_f$, we need two additional subsets:
\[
    \cA_{I}^+ \coloneqq \big\{ 
        \alpha_i \in \cA \mid \alpha_i \sim J \text{ and } J \leq I 
    \big\},
\]
\[
    \cA_{I}^- \coloneqq \big\{ 
        \alpha_i \in \cA \mid \alpha_i \sim J \text{ and } J < I 
    \big\}.
\]
Then, the value of $\cM_f(I,J)$ is calculated using the restriction of the linear map $f_d$, where $d$ is the right endpoint of $J$, to the subspaces:
    \[
        F_{IJd} \coloneqq
        \left(
        \begin{array}{c|cc}
         &
        \cA^-_{Id} &
        \cA_{Id} 
        \\
        \hline
        \cB_{Jd} &
        *  &
        * 
        \\
        \cB_d \setminus \hcB^+_{Jd} &
        *  &
        * 
        \end{array}
        \right).
    \]
Then, we proceed with a column reduction of the matrix. 
$\cM_f(I,J)$ is given by the number of pivots with row in $\cB_{Jd}$ and column in $\cA_{Id}$. 
\smallskip

In general, $\matching{f}$ and $\cM(f)$ are not equal.
We provide now an example of this, leaving the calculations to the reader.
\begin{example}
	Consider again Example~\ref{ex:complex_indecomposable}.
	and recall that $V\cong k_{[1,4]}\oplus k_{[2,3]}$ and $U\cong k_{[1,3]}\oplus k_{[2,2]}$. 
	Using the matrix method for calculating $\cM_f$ explained above, we have that $\cM_f([1,4], [1,3])= 1$, $\cM_f([2,3], [1,3])= 1$ and 
	$\cM_f(I, J)= 0$ otherwise.  
\end{example}

However, both operators are the same for interval decomposable maps.

\begin{proposition}
 If $f$ is an interval decomposable persistence map, then  $\matching{f} = \cM(f)$.
\end{proposition}
\begin{proof}
    Both operators are additive, and if $f_{IJ} \colon k_{I} \to k_{J}$ is a non-null map, their only non-null value is $\matching{f}(I,J) = \cM_f(I,J) = 1$.
\end{proof}
In \citep[Theorem~5.3]{matrix_oxford}, it was proven that any $f \colon V \to U$ for which $V$ and $U$ does not have nested intervals is necessarily interval decomposable.
Hence, we have the following corollary.
\begin{corollary}
     If $f \colon V \to U$ is a persistence map such that there are no nested intervals in $\II(V)$ or $\II(U)$, then $\matching{f} = \cM_f$.   
\end{corollary}
Finally, we provide an example of how to combine the outputs of $\matching{f}$ and $\cM_f$ to prove that a module is indecomposable.
\begin{example}
	Consider again Example~\ref{ex:complex_indecomposable}
	and recall that $V\cong k_{[1,4]}\oplus k_{[2,3]}$ and $U\cong k_{[1,3]}\oplus k_{[2,2]}$. 
	We show now that $f$ is indecomposable by using the outputs of $\cM_f$ and $\matching{f}$. Suppose that 
    \[ 
        f \cong f' \oplus f'' \colon V' \oplus V'' \to U' \oplus U''.
    \] 
	Without loss of generality, we assume that $V'\neq 0$ and that $k_{[1,4]}$ is a summand in the decomposition of $V'$ and not in $V''$.
	Since $\cM_f$ is additive, from $\cM_f([1,4], [1,3])= 1$ we get that $k_{[1,3]}$ must be a summand of $U'$.
	In addition, $\cM_f([2,3], [1,3])= 1$ implies that $k_{[2,3]}$ is
	a summand of $V'$, and necessarily $V'' = 0$.
	Since $\matching{f}([2,3], [2,2])=1$ and $\matching{f}$ is also additive, $k_{[2,2]}$ must also be a summand of $U'$.
	Hence, $U''= 0$ and $f'':V''\rightarrow U''$ is trivial, so $f = f'$ and it is indecomposable.
\end{example}

\end{document}